\documentclass{article}

\usepackage{geometry}
\usepackage{cmbright} 
\usepackage[T1]{fontenc} 
\usepackage[utf8]{inputenc} 
\usepackage{amsmath}
\usepackage{amsthm}
\usepackage{amsfonts}
\usepackage{graphicx}
\usepackage{algorithm}
\usepackage{algorithmic}
\usepackage{bm}

\usepackage{biblatex}
\graphicspath{{pictures/}}

\newtheorem{remark}{Remark}

\DeclareMathAlphabet\mathbfcal{OMS}{cmsy}{b}{n}

\newcommand{\bfA}{ {\bf A} }
\newcommand{\bfu}{ {\bf u} }

\newcommand{\bfj}{  {\boldsymbol{j}  } }
\newcommand{\ibfj}{ {\boldsymbol{j}  } }
\newcommand{\bfxi}{ {\bm \xi} }
\newcommand{\bff}{ {\bf f} }
\newcommand{\bfv}{ {\bf v} }
\newcommand{\bfH}{ {\bf H} }
\newcommand{\bfzero}{ {\bm 0} }

\newcommand{\R}{ \mathbb R }
\newcommand{\E}{ \mathbb E }
\newcommand{\diff}{ \; \text{d} }
\newcommand{\var}{ \text{Var} }
\newcommand{\minmod}{ \text{minmod} }

\newcommand{\flux}[4]{ {\mathbfcal F}_{ {#1},{#2} }^{ {#3} {#4}} }
\newcommand{\fluxE}[4]{ \bar{\mathbfcal F}_{ {#1},{#2} }^{ {#3} {#4}} }

\newcommand{\matL}{ {\bf L} }
\newcommand{\matR}{ {\bf R} }
\newcommand{\matI}{ {\bf I} }
\newcommand{\rmax}{ {r_\text{max}} }
\newcommand{\Call}[2]{\text{#1}\,(#2)}

\title{High-dimensional stochastic finite volumes using the tensor train format}
\date{}
\author{Juliette Dubois\footnote{1}, Michael Herty\footnote{2}, Siegfried M\"uller\footnote{2} }

\author{
\renewcommand{\thefootnote}{\arabic{footnote}} 
Juliette Dubois\footnotemark{} 
\renewcommand{\thefootnote}{\fnsymbol{footnote}} 
\footnotemark[1]{}
\renewcommand{\thefootnote}{\arabic{footnote}},
Michael Herty\footnotemark[2]{} ,
\renewcommand{\thefootnote}{\arabic{footnote}} 
Siegfried M\"uller\footnotemark[2]{}
}

\begin{document}
\maketitle

\footnotetext[1]{Institut für Mathematik, Technische Universität Berlin, Straße des 17. Juni 136, 10623 Berlin, Germany}
\renewcommand{\thefootnote}{\fnsymbol{footnote}} 
\footnotetext[1]{Corresponding author: dubois@math.tu-berlin.de} 
\renewcommand{\thefootnote}{\arabic{footnote}} 

\footnotetext[2]{Institut für Geometrie und Praktische Mathematik, RWTH Aachen University, Templergraben 55, D-52056 Aachen}

\paragraph{Abstract. }
A method for the uncertainty quantification of nonlinear hyperbolic conservation laws with many uncertain parameters is presented. 
The method combines  stochastic finite volume methods and tensor trains in a novel way: 
the dimensions of physical space and time  are kept as full tensors, while all stochastic dimensions are compressed together into a tensor train. 
The resulting hybrid format has one tensor train for each spatial cell and each time step. \\
The MUSCL scheme is adapted to the proposed hybrid format, and its feasibility is demonstrated through several test cases.  
For the scalar Burgers’ equation, we conduct a convergence study and compare the results with those obtained using the full tensor train format with three stochastic parameters. The equation is then solved for an increasing number of stochastic dimensions.  
For systems of conservation laws, we focus on the Euler equations. A parameter study and a comparison with the full tensor train format are carried out for the Sod shock tube problem. As a more complex application, we investigate the Shu–Osher problem, which involves intricate wave interactions.  
The presented method opens new avenues for integrating uncertainty quantification with established numerical schemes for hyperbolic conservation laws.


\paragraph{MSCcodes.} 65M08, 65M75, 35R60

\section{Introduction}
This work addresses uncertainty quantification for nonlinear hyperbolic conservation laws with multiple parameters.  
A variety of methods have been developed to analyze how uncertainties influence partial differential equations.  
Two principal approaches are Monte Carlo-type methods and methods based on orthogonal expansions, such as the generalized polynomial chaos (gPC) framework.  
The gPC approach can be implemented either intrusively, for example through stochastic Galerkin methods, or non-intrusively, such as stochastic collocation.  
Monte Carlo methods are non-intrusive and their convergence rate is independent of the problem dimension; however, they converge very slowly.  
In contrast, gPC methods achieve exponential convergence when the solution depends smoothly on the stochastic parameters.  
This property makes gPC particularly effective for elliptic and parabolic equations, but such smooth dependence is typically not observed for solutions to stochastic hyperbolic partial differential equations.

In this context, the stochastic finite volume method (SFV) was introduced in \cite{barth_propagation_2012} to study conservation laws with uncertainties. 
The SFV is a deterministic formulation of the equations that keeps some properties of the original hyperbolic problem, such as well-posedness \cite{mishra_sparse_2012}. 
However, since a new dimension is added for each uncertain parameter, solving the equations numerically becomes infeasible. 

When methods originally developed for low-dimensional settings are extended to high dimensions, the number of parameters grows exponentially. This phenomenon is known as the curse of dimensionality.  
Low-rank tensor formats offer a potential remedy by keeping storage requirements and computational operations manageable. They can be viewed as a generalization of the well-known low-rank matrix decomposition to higher-order tensors.  
Several tensor formats have been introduced, including the canonical, Tucker, and hierarchical formats \cite{bachmayr_low-rank_2023}. In this work, we focus on one hierarchical format: the tensor train format.  Tensor trains have been successfully applied to elliptic, parabolic, and linear hyperbolic problems \cite{oseledets_tensor-train_2011, zhong_fast_2018}.  
Nonlinear hyperbolic conservation laws, however, present additional challenges, due to the development of shocks in the solution.   In \cite{walton_tensor-train_2025}, the authors propose combining SFV with tensor trains and report promising results.

In this paper, we propose a hybrid formulation where --- in contrast to  \cite{walton_tensor-train_2025} ---  time and physical space are kept in the full format while the stochastic space is compressed in the tensor train format.
Section \ref{sec:SFV} introduces the problem and presents the principle of the stochastic finite volume method.
In Section \ref{sec:TT} we recall the main ideas of the tensor train format.
The new hybrid format is presented in Section \ref{sec:hybrid}.
We describe a MUSCL-type algorithm adapted to the hybrid format.  
The presentation emphasizes the minimal modifications required to incorporate tensor trains compared to the classical scheme.
In Section \ref{sec:Numerics}, we show numerical experiments that prove the feasibility and efficiency of our approach.

\section{The stochastic finite volume method} \label{sec:SFV}
We are interested in conservation laws with uncertain initial data. 
For ease of presentation, the problem is presented in the spatially one-dimensional case. We emphasize that the method is the same for higher dimensions. 
Let $T>0$, $\Omega_x\subset\R$, $\Omega_\bfxi\subset\R^m$, 
and let $\bfu:(0,T)\times\Omega_x\times\Omega_\bfxi\to\mathbb{R}^p$ be the solution to 
\begin{align}  
    \label{eq:PDE}
    \frac{\partial \bfu }{\partial t}(t,x; \omega_\bfxi)
    + \frac{\partial \bff}{\partial x} (\bfu (t,x; \omega_\bfxi))
   & = \bfzero, \quad x \in \Omega_x, \ \omega_\bfxi \in \Omega_\bfxi, \ t \in (0,T),
    \\
    \label{eq:PDE_IC}
    \bfu (0, x; \omega_\bfxi) & = \bfu_0(x; \omega_\bfxi)
    , \quad x \in \Omega_x, \ \omega_\bfxi \in \Omega_\bfxi.
\end{align}
Here, $\bff = (f_1, \dots f_p)$ is the flux field. 
The random input $\omega_\bfxi$ is parametrized by a random variable $\bfxi:\Omega\to\Omega_\bfxi$ defined on the probability space 
$(\Omega, \mathcal F, \mathbb P)$.
We assume that there exists a probability density 
$p:\R^m\to[0,\infty)$ such that the expectation and the variance of $\bfu$ can be expressed as 
\begin{equation}
\E[\bfu(t,x)] = \int_{\Omega_\bfxi} \bfu(t,x ; \omega_\bfxi) p (\bfxi) \diff \bfxi, \quad
\var [\bfu(t,x)] = \E[(\bfu(t,x)-\E[\bfu(t,x)])^2].  
\end{equation}

\subsection{General description of the method}
The idea is to consider the space $\Omega_\bfxi\subset\R^m$ as new ``spatial directions'' of the problem \cite{jin_stochastic_2017}. 
The stochastic problem \eqref{eq:PDE}-\eqref{eq:PDE_IC} is reformulated as a deterministic problem for the unknown 
$\bfu(t,x,\bfxi)$. 

The physical space $\Omega_x$ and the parametrized probability space $\Omega_\bfxi$
are discretized by Cartesian grids denoted by $\mathcal{C}_x$ and $\mathcal{C}_\bfxi$, respectively. 
The time interval $[0,T]$ is discretized with a time step $\Delta t$, 
\[
\mathcal C_x   = \cup_i K_x^i, 
\quad
\mathcal C_\bfxi = \cup_\ibfj K_\bfxi^\ibfj, 
\quad
[0,t]= \cup_k [t_k, t_k + \Delta t], 
\]
where $\bfj=(j_1, \dots j_m)$ is a multi-index.
We introduce a constant mesh size in the spatial direction $\Delta x$,
and in each stochastic dimension $\Delta \xi_1,\dots,\Delta \xi_m$.  
The cells $K_x^i$ and $K_\bfxi^\ibfj$ are defined as 
\begin{align}
    K_x^i = [x_{i-1/2}, x_{i+1/2} ], 
    \quad x_{i \pm 1/2} = x_{i} \pm \frac{\Delta x}{2}, 
    \\
    K_\bfxi^\ibfj = \Pi_{\ell = 1}^m [\xi_{j_\ell-1/2}, \xi_{j_\ell + 1/2} ],
    \quad 
    \xi_{j_\ell\pm 1/2} = \xi_{j_\ell } \pm \frac{\Delta \xi_\ell}{2} .
\end{align}
Let $N_x$ and $N_{\xi_1},\dots,N_{\xi_m}$ denote the number of 1D cells in each spatial and stochastic direction, respectively. 

We introduce the cell average in space
\begin{equation}
    \bfu_{i}^n(\bfxi) =
    \frac{1}{|K_x^i|}
    \int_{K_x^i}  \bfu (t^n,x,\bfxi) \diff x,
\end{equation}
and the cell average in space and expectation over a cell $\bfj$
\begin{equation}
    \bar \bfu_{i,\ibfj}^n =
    \frac{1}{|K_\bfxi^\ibfj|}
    \E_\ibfj [\bfu_{i,\ibfj}^n]
    =
    \frac{1}{|K_x^i||K_\bfxi^\ibfj|}
    \int_{K_\bfxi^\ibfj}  \int_{K_x^i}  \bfu (t^n ,x,\bfxi)
    \diff x \, p(\bfxi) \diff \bfxi.
\end{equation}
Here, the cell volumes are 
\begin{equation}
    |K_x^i| = \int_{K_x^i} \diff x = \Delta x, 
    \qquad 
    |K_\bfxi^\ibfj| = \int_{K_\bfxi^\ibfj} p(\bfxi) \diff \bfxi.
\end{equation}
\\

Integrating Equation \eqref{eq:PDE} in space and computing the expectation in each cell yields   
\begin{equation} \label{eq:PDE_integral}
    \int_{K_\bfxi^\bfj} \int_{K_x^i}
    \frac{\partial \bfu }{\partial t}(t^n ,x,\bfxi)
     \diff x \, p(\bfxi) \diff \bfxi 
    + 
    \int_{K_\bfxi^\bfj} \int_{K_x^i} 
    \frac{\partial \bff}{\partial x} (\bfu (t^n ,x,\bfxi))
     \diff x  \, p(\bfxi) \diff \bfxi 
    = \bfzero. 
\end{equation}
After integration by parts, Equation \eqref{eq:PDE_integral} reads 
\begin{equation}\label{eq:PDE_integral-2}
    \Delta x |K_\bfxi^\ibfj| \frac{d \bar \bfu_{i,\ibfj}^n}{d t}
    + \E_\ibfj [ {\bf f}(\bfu (t,x_{i+1/2},\bfxi) ) - {\bf f}(\bfu (t,x_{i-1/2},\bfxi)) ]
    = \bfzero.
\end{equation}

The exact flux $\bff$ is replaced by a numerical flux $\mathbfcal F$: 
the value of the solution at the interface $x_{i+1/2}$ is 
approximated using the values of the averages in the surrounding cells 
$\bar \bfu_{i-1,\ibfj}^n , \bar \bfu_{i,\ibfj}^n , \bar \bfu_{i+1,\ibfj}^n $. 
Let $\flux{i}{\ibfj}{n}{+}$ and $\flux{i}{\ibfj}{n}{-}$ denote the numerical flux at the interface 
$x_{i+1/2}$ and $x_{i-1/2}$, respectively. 
The approximation of their expectation on each cell is denoted by $\fluxE{i}{\ibfj}{n}{+}$ and $\fluxE{i}{\ibfj}{n}{-}$. Then Equation \eqref{eq:PDE_integral-2} is approximated by
\begin{equation} \label{eq:PDE_numFlux}
    \Delta x |K_\bfxi^\ibfj| \frac{d \bar \bfu_{i,\ibfj}^n }{d t}
    + \fluxE{i}{\ibfj}{n}{+} - \fluxE{i}{\ibfj}{n}{-}
    = 0, 
    \qquad \fluxE{i}{\ibfj}{n}{+}
    = \E_\ibfj [ \mathcal F^+(\bfu_{i,\ibfj}^n ) ] , 
    \quad \fluxE{i}{\ibfj}{n}{-}
    = \E_\ibfj [ \mathcal F^-(\bfu_{i,\ibfj}^n ) ].
\end{equation}

The integrals are approximated using the midpoint rule. 
Then the approximated average over a spatial cell reads
\begin{equation} \label{eq:quad_space}
\bar \bfu_{i,\ibfj}^n  \approx 
\bfu(t^n , x_i, \bfxi_\ibfj) 
\end{equation}
By a slight abuse of notation, let $\bar \bfu_{i,\ibfj}^n$ and $\fluxE{i}{\ibfj}{n}{\pm}$ denote  the expectation of the solution and the numerical fluxes, respectively, both approximated by quadrature. The midpoint rule is used to compute the expectation $\mathbb{E}$ with respect to $p$, and, therefore, quantities are evaluated at the point $\bfxi_\ibfj$. 
We introduce the approximated volume on the stochastic cells $K_\ibfj=\Delta\xi_1\dots\Delta\xi_m\,p(\xi_\ibfj)$.
Then the semi-discrete equation reads 
\begin{equation} \label{eq:PDE_semiDiscrete}
    K_\ibfj \Delta x \frac{d \bar \bfu_{i,\ibfj}^n }{d t}
    + K_\ibfj \fluxE{i}{\ibfj}{n}{+}
    - K_\ibfj \fluxE{i}{\ibfj}{n}{-}
    = \bfzero.
\end{equation}
After simplifications ($K_\ibfj\neq 0$), the semi-discrete scheme is determined by
\begin{equation}
    \frac{d \bar \bfu_{i,\ibfj}^n }{d t}
    + \frac{1}{\Delta x} \left(
        \fluxE{i}{\ibfj}{n}{+} - \fluxE{i}{\ibfj}{n}{-}
    \right)
    = \bfzero.
\end{equation}

\subsection{Numerical flux}
The numerical flux is computed for a MUSCL scheme \cite{kurganov_new_2000}.
The expectation of the numerical flux at the interface $i+1/2$ is defined by
\begin{equation} \label{eq:Full_NumFlux}
\bar{\bfH}_{i+1/2,\ibfj}^n  = 
\frac{ \bff \left( \bfu_{i+1/2,\ibfj}^{n +} \right) 
     + \bff \left( \bfu_{i+1/2,\ibfj}^{n -} \right) }{2}
- \frac{a_{i+1/2,\ibfj}^n}{2} \left( \bfu_{i+1/2,\ibfj}^{n +}  - \bfu_{i+1/2,\ibfj}^{n -}  \right),
\end{equation}
where intermediate values $\bfu_{i+1/2,\ibfj}^{n +}, \bfu_{i+1/2,\ibfj}^{n -}$ are linear reconstructions of the solution using the approximate derivative,
and $a_{i+1/2,\ibfj}^n$ is the local speed. 

The intermediate values are computed at the interface between two cells in space: 
\begin{align} \label{eq:Full_reconstruction}
    \bfu_{i+1/2,\ibfj}^{n +} &= \bar \bfu_{i+1,\ibfj}^n  - \frac{\Delta x}{2} (\bfu_x)_{i+1,\ibfj}^n  
    \\
    \bfu_{i+1/2,\ibfj}^{n -} &= \bar \bfu_{i,\ibfj}^n  + \frac{\Delta x}{2} (\bfu_x)_{i,\ibfj}^n  
\end{align}
The approximation of the derivative $ (\bfu_x)_{i,\ibfj}^n $ is obtained using a classical $\minmod$ limiter,
\begin{equation} \label{eq:Full_deriv}
    (\bfu_x)_{i,\ibfj}^n  = \minmod \left(
        \frac{\bar{\bfu}_{i,\ibfj}^n - \bar{\bfu}_{i-1,\ibfj}^n }{\Delta x}, \frac{\bar{\bfu}_{i+1,\ibfj}^n - \bar{\bfu}_{i,\ibfj}^n }{\Delta x}
    \right). 
\end{equation}
The local speed is given by 
\begin{equation}
    a_{i+1/2,\ibfj}^n  = \max\left\{
        \rho(\bff'(\bfu_{i+1/2,\ibfj}^{n +})), \rho(\bff'(\bfu_{i+1/2,\ibfj}^{n +}))
    \right\},
\end{equation}
where $\rho(\bfA)$ denotes the spectral radius of the matrix $\bfA$. 
For a matrix $\bfA\in\R^{N\times N}$ with eigenvalues $(\lambda_i)_{1\leq i \leq N}$, we recall that the spectral radius is defined by 
$\rho(\bfA) = \max_{1\leq i \leq N} |\lambda_i|$.

\subsection{Fully-discrete formulation and computational cost}
We conclude the presentation of the stochastic finite volume method with the time discretization. 
A simple forward Euler scheme yields
\begin{equation} \label{eq:Full_discrete}
    \bar \bfu_{i,\ibfj}^{n +1} 
    = \bar \bfu_{i,\ibfj}^n  
    - \frac{\Delta t}{\Delta x} (\bar \bfH_{i+1/2,\ibfj}^n  - \bar \bfH_{i-1/2,\ibfj}^n ).
\end{equation}
Note that in the case of a forward Euler discretization in time, the cell average in space is computed only for the initial condition. 
The cell average of the next time steps is directly given by \eqref{eq:Full_discrete}.

The stochastic finite volume method provides a complete description of the dependence of the solution on  the random variables. 
However, the problem becomes computationally intractable with an increasing  dimension of the stochastic space $\Omega_\bfxi$. 
Indeed, in the general case 
$m\geq 1$ with
$N_{\xi_1},\dots,N_{\xi_m}$ stochastic cells and $N_t$ time steps, 
the discrete solution $\bar \bfu_{i,\ibfj}^n $ is a tensor whose number of entries grows exponentially with the number of dimensions: 
\begin{equation}
    \left(\bar \bfu_{i,\ibfj}^n \right)_{
        \begin{subarray}{l}n=1,\dots,N_t \\
            i=1,\dots,N_x \\
            \ibfj = (1,\dots,1),\dots,(N_{\xi_1},\dots,N_{\xi_m})
            \end{subarray}
        }
    \in \R^{p N_t} \times 
    \R^{p N_{x_1}}   \dots \times \R^{p N_{x_d}} \times 
    \R^{p N_{\xi_1}} \dots \times \R^{p N_{\xi_m}}.
\end{equation}
For the  problem considered, 
$m$  becomes very large, such that the tensor $\bar \bfu_{i,\ibfj}^n$ is too large for both computation and storage. 
To circumvent this problem, we propose to use a low-rank tensor approximation. 
Following \cite{walton_tensor-train_2025} we consider  the tensor train format to represent the low-rank structure. The tensor trains are briefly described in the next section. 

\section{The tensor train format} \label{sec:TT}
Various low-rank formats for tensors have been proposed \cite{bachmayr_low-rank_2023}. 
We focus on the tensor train format, that has been introduced by \cite{oseledets_tensor-train_2011}.  

\subsection{Definition of a tensor train}
For $m \in \mathbb N$, a $m$th-order tensor $T~\in~\R^{N_1 \times \dots \times N_m}$ is a tensor train if it can be written as a sum of products of third-order tensors,
\begin{equation} \label{eq:Def_TT}
    T(i_1,\dots,i_m) 
    = \sum_{\alpha_1, \dots, \alpha_m =1 }^{r_1, \dots, r_m}
    G_1(i_1,\alpha_1) G_2(\alpha_1, i_2, \alpha_2) \cdots G_m(\alpha_m, i_m).
\end{equation}
The third-order tensors $G_\ell \in \R^{r_{\ell-1} \times N_\ell \times r_\ell}$ are called the cores of the tensor train, and the values $r_1,\dots,r_m$ are called the ranks.
For the first and the last core, the ranks are set to $r_0=r_{m+1}=1$, respectively.
The tensor train rank (TT-rank) is the largest rank of the tensor train. We denote it by $r=\max_{\ell=1 \dots m}(r_\ell)$. 

The tensor train format has two properties that make it well suited for our objective: it breaks the curse of dimensionality, and any tensor has a quasi-optimal approximation in the tensor-train format. 
The first property is the primary  goal of using low-rank formats. 
For the $m$th-order tensor with $N$ entries in each dimension, a total of $N^m$ entries is required to describe the tensor. 
Adding new dimensions to the tensor exponentially increases  the number of entries. 
However, for a tensor train, only $\mathcal O((m-2)Nr^2+2Nr)$ entries are necessary. Therefore, the number of entries is only polynomial in the number of dimensions. 
The second property fulfills a practical requirement: not all tensors have a low-rank representation, thus it is  necessary to  approximate a given tensor by a tensor in the low-rank format. 
It was proven in  \cite{oseledets_tensor-train_2011} that the problem of finding a quasi-optimal tensor train approximation for a given rank is a well-posed problem. 
The proof is constructive, and  an algorithm to construct the quasi-optimal approximation is provided. 
We stress here that such algorithm is not available for all low-rank formats. For example, the canonical decomposition has a higher compression rate, but its approximation problem is ill-posed \cite{de_silva_tensor_2008}. We therefore will use the tensor train format in the following. 

\subsection{Operations on tensor trains}

Tensor trains allow for a range of algorithms for many common tensor operations. 
We list here the ones that will be used in the rest of the paper. 
Algebraic operations on tensors, such as addition, element-wise multiplication, and tensor contraction can be performed directly in the tensor train format by applying the appropriate operations on the cores. 
However, the resulting tensors have larger ranks after the operation.
Therefore, a rounding algorithm was introduced in \cite{oseledets_tensor-train_2011}. 
Rounding a tensor train yields a tensor train approximation with lower rank. 
After each algebraic operation, the result is rounded to keep the rank low.

Applying a non-algebraic function to a tensor train is also possible using  the cross approximation algorithm. The cross approximation was first introduced in \cite{oseledets_tt-cross_2010}, and a version with adaptive ranks was developed in \cite{savostyanov_fast_2011}.
With cross approximation, the non-algebraic function is evaluated at chosen entries of the tensor, and a tensor train approximation is deduced from those entries.
Error estimates for the cross approximation algorithm have been developed  for element wise errors, cf.~\cite{savostyanov_quasioptimality_2014,osinsky_tensor_2019}, and   for errors in the Frobenius norm, cf.~\cite{qin_error_2022}. 
The main result of those papers is that the approximation error does not grow exponentially with the number of dimensions. 


\section{A hybrid format for the stochastic finite volume method} \label{sec:hybrid}
We propose a method for solving hyperbolic conservation laws with a large number of uncertain parameters.  
Our approach exploits low-rank approximation while preserving the structure of deterministic solvers.  
Specifically, the tensor train format is employed only in the stochastic space, while the full format is retained in the spatial and temporal domains.  
We deliberately keep the full format in time, since the solution is not expected to admit a low-rank representation with respect to this variable. In nonlinear hyperbolic problems, the solution can deviate substantially from the initial condition as time evolves, due to gradient steepening and the subsequent development of discontinuities.
Comparisons between low-rank approximations and full formats have been reported, for example in \cite{zhong_fast_2018}, for linear and quasi-linear hyperbolic problems. The authors observed that retaining the full format in time yields higher accuracy than its low-rank counterpart, and this effect can be expected to be even more pronounced in nonlinear problems.
Moreover, maintaining the full format in both space and time facilitates the adaptation of well-established deterministic solvers, since most of the implementation follows classical algorithms. Scalar operations are simply replaced by operations on tensor trains, as described in more detail below.

The low-rank format is used only in the  stochastic directions. 
However, there is no theoretical guarantee that the solution to an initial value problem with low-rank structure will itself retain a low-rank representation. In fact, previous studies \cite{pares_higher-dimensional_2024, Kolb:2024} have observed the development of discontinuities in the stochastic dimension, suggesting that a locally full-rank approximation may be required.
Despite this limitation, other results \cite{walton_tensor-train_2025} on related problems provide experimental evidence that low-rank formats remain a viable tool for approximating solutions. The method introduced in that work compresses both the physical and stochastic spaces into a single tensor train. This representation, referred to hereafter as the full-TT format, will serve as a basis for comparison.
\par 
We now discuss the proposed method. The method is stated for a scalar equation in a single spatial dimension and using the same discretization in each stochastic dimension. This is done for simplicity of the presentation only. 
We define 
$\Delta\xi:= \Delta\xi_1=\Delta\xi_2=\dots=\Delta\xi_m$, and 
$N_\xi:= N_{\xi_1}=N_{\xi_2}=\dots=N_{\xi_m}$. 
A tensor train is used to approximate the solution for each spatial cell $i$, hence, the computational variable is an array of tensor trains:
\begin{equation} 
\bar{U}^n=(\bar{U}_1^n,\dots,\bar{U}_{N_x}^n)^\top, 
\quad \text{ with } \quad
\bar U_i^n  \approx  \left( \bar u_{i,\ibfj}^n  \right)_{\ibfj \in \{1,\dots,N_\xi\}^m}.
\end{equation}
This hybrid-TT format can be described by adding a new dimension to each core,
\begin{equation}
    \bar{U}^n_i (\ibfj) 
    = \sum_{\alpha_1, \dots, \alpha_m =1 }^{r_1, \dots, r_m}
    G_1(i,j_1,\alpha_1) G_2(i,\alpha_1, j_2, \alpha_2) \cdots G_m(i,\alpha_m, j_m), 
\end{equation}
however, we will rather consider the hybrid-TT format as an array of tensor trains in the following. 
For comparison, we state the full-TT format. A single tensor train with a core for the spatial variable is used:
\begin{equation} \label{eq:FullyTT_def}
    \bar{U}^n(i,\ibfj) 
    = \sum_{\alpha_1, \dots, \alpha_m =1 }^{r_1, \dots, r_m}
    G_0(i,\alpha_0) G_1(\alpha_0,j_1,\alpha_1) G_2(\alpha_1, j_2, \alpha_2) \cdots G_m(\alpha_m, j_m).
\end{equation}

The approach is exemplified for the MUSCL finite--volume scheme and a  first--order explicit time integration. The derivation can be easily generalized and it may be applied to any other finite-volume scheme. We also recall the equations for the evaluation  of the expectation and the variance of the solution computed by the hybrid and full-TT schemes, respectively. 

\subsection{The MUSCL scheme in the hybrid format}
The scheme of Section \ref{sec:SFV} is adapted to the hybrid format. 
The reconstructed states and the expectation of the numerical flux in the hybrid format are  denoted by $U^{n\pm}=(U_{3/2}^{n\pm},\dots,U_{N_x+1/2}^{n\pm})^\top$, 
and $\bar{H}^n=(\bar{H}^n_{1/2},\dots,\bar{H}^n_{N_x+1/2})^\top$, respectively. 
We assume that for a fixed cell $i$ and time step $n$ the tensor train approximations
$\bar U^n _i$ and  $U_{i+1/2}^{n +}$,  $U_{i+1/2}^{n -}$ are given. 
Then, the MUSCL scheme in hybrid format reads:
\begin{align}
    \label{eq:FV-TT1}
    \bar U_i^{n +1} 
    & = \bar U_i^n  - \frac{\Delta t}{\Delta x}  (\bar H_{i+1/2}^n  - \bar H_{i-1/2}^n ),
    \\[5pt]
    \bar H_{i+1/2}^n  
    & = \frac{ f\left( U_{i+1/2}^{n +} \right) + f\left( U_{i+1/2}^{n -} \right) }{2}
    - \frac{a_{i+1/2}}{2} \left(  U_{i+1/2}^{n +}  -  U_{i+1/2}^{n -}  \right),
    \\[5pt] 
    a_{i+1/2}^n
    & = \max \left(
        |f'(U_{i+1/2}^{n +})|, |f'(U_{i+1/2}^{n -})|
    \right), 
    \\
    U_{i+1/2}^{n +} 
    & = \bar U_{i+1}^n  - \frac{\Delta x}{2} (U_x)_{i+1}^n, 
    \\
    U_{i+1/2}^{n -} 
    & = \bar U_{i}^n  + \frac{\Delta x}{2} (U_x)_{i}^n, 
    \\
    (U_x)_{i}^n &= \minmod\left(
        \frac{\bar U_i^n - \bar U_{i-1}^n}{\Delta x},
        \frac{\bar U_{i+1}^n - \bar U_i^n}{\Delta x}
    \right) .
    \label{eq:FV-TT5}
\end{align}
The equations \eqref{eq:FV-TT1}-\eqref{eq:FV-TT5} are similar to a classical MUSCL scheme without stochastic variables.  The only difference is that the objects 
$\bar U_i^n , U_{i+1/2}^{n \pm}, (U_x)_i^n $ and $\bar H_{i+1/2}^n $ 
are scalars in the classical case and tensor trains in our hybrid formulation. 
Hence, all operations on these objects must be adapted to the tensor train format. 
As outlined in Section \ref{sec:TT}, linear and polynomial functions may be applied directly to tensor trains, whereas the cross-interpolation algorithm must be used for non-polynomial functions. 
In our case, the operations $(T_1, T_2) \mapsto \minmod(T_1, T_2)$ and 
the function $(T_1, T_2) \mapsto \max(|T_1|, |T_2|)$ are non-polynomial in the variables $T_1, T_2$. 
Depending on the problem at hand, this can also be the case for the initial condition $u_0$, the flux $f$ and its derivative $f'$. 

The first step of the algorithm is to compute the cell average of the initial condition in the hybrid format $\bar U^0$ using \eqref{eq:quad_space}. 
In the  scalar case, the maximum principle $|u(t,x,\bfxi)| \leq \max_{x,\bfxi}|u_0(x,\bfxi)|$ holds and we define  the fixed time step
\begin{equation} \label{eq:CFL}
    \Delta t < \frac{1}{2} \frac{\Delta x}{\max|f'(\bar U^0)|}. 
\end{equation}
We refer to Section \ref{sec:Euler} for  systems of conservation laws.
The solution is then propagated in time using the scheme \eqref{eq:FV-TT1}-\eqref{eq:FV-TT5}. The algorithm  is summarized in the Appendix \ref{sec:Algo}.

\subsection{Computing expectation and variance} \label{sec:Exp_var}
The expectation and variance of a solution in the hybrid format $\bar U$ is computed as follows. We describe the computation for two stochastic variables $\bfxi=(\xi_1, \xi_2)$, the case $m>2$ being similar. 
The expectation and variance of the solution $u$ at time $t$ and point $x$ are 
\begin{align} 
\label{eq:exp}
    \E[u(t,x)] &= \int_{\Omega_\bfxi} u(t,x,\xi_1,\xi_2) p_{\xi_1,\xi_2}(\xi_1,\xi_2) \diff \xi_1 \diff \xi_2, 
    \\
\label{eq:var}
    \var [u(t,x)] &= \E[(u(t,x)-\E[u(t,x)])^2] = \E[u(t,x)^2] - \E[u(t,x)]^2.  
\end{align}
We make the additional assumption that the random variables are independent and denote by $p_1, p_2$ the probability distributions associated with $\xi_1$ and $\xi_2$, respectively. 
The integral on the whole space $\Omega_\bfxi$ is approximated as a sum 
\begin{equation}
    \E[u_i^n] \approx \sum_{j_1=1}^{N_\xi} \sum_{j_2=1}^{N_\xi} 
    \bar u_{i,j_1,j_2}^n p_1(\xi_{j_1}) p_2(\xi_{j_2}) (\Delta \xi)^2.
\end{equation}
This operation can be seen as a contraction of the tensor 
$\bar u_{i,j_1,j_2}^n$ with the tensor $p(\xi_{j_1})p(\xi_{j_2})$. 
Moreover, the tensor $p(\xi_{j_1})p(\xi_{j_2})$ is a tensor train with rank one and cores $p(\xi_{j_1})$ and $p(\xi_{j_2})$. 
The expectation is thus computed in the tensor train format.
The variance is computed using Equation \eqref{eq:var}.  We first compute the pointwise square of the solution  in the hybrid format by an operation on tensor trains, then its expectation is obtained as above.

\subsection{The full-TT format for the MUSCL scheme} 
We compare the proposed approach with the full-TT approach. To this end, the method introduced in \cite{walton_tensor-train_2025} is stated for the  MUSCL scheme. 
The core idea is to rewrite the scheme \eqref{eq:Full_NumFlux}-\eqref{eq:Full_discrete} in terms of global matrices. Then, every computational operation is applied to the full-TT approximation of the solution, using either TT algebraic operations or cross-approximation. 

We illustrate the idea by describing the reconstruction step using global matrices. 
For simplicity, the operation is first described in the scalar case and  one spatial dimension. Let $\bar u^n = (u_{1}^n, \dots u_{N_x}^n)$ be the solution vector. 
We introduce the identity matrix $\matI \in \mathbb R^{N_x \times N_x}$ and the matrices $\matL, \matR \in \mathbb R^{N_x \times N_x}$ defined by
\begin{align*}
    &\matL = \begin{pmatrix}
        0      & 0      & 0      & \cdots & 0      \\
        1      & 0      & 0      & \cdots & 0      \\
        0      & 1      & 0      & \cdots & 0      \\
        \vdots & \vdots & \vdots & \vdots & \vdots \\
        0      & 0      & \cdots & 0      & 0      \\
        0      & 0      & \cdots & 1      & 0      \\
    \end{pmatrix}  \text{ and }
    \matR = \begin{pmatrix}
        0      & 1      & 0      & \cdots & 0      \\
        0      & 0      & 1      & \cdots & 0      \\
        0      & 0      & 0      & \cdots & 0      \\
        \vdots & \vdots & \vdots & \vdots & \vdots \\
        0      & 0      & \cdots & 0      & 1      \\
        0      & 0      & \cdots & 0      & 0      \\
    \end{pmatrix} .
\end{align*}
The vector of approximate derivatives defined in \eqref{eq:Full_deriv} is obtained as a matrix--vector multiplication with those matrices, i.e.,  
\begin{equation} \label{eq:fullTT_deriv}
    (u_x)^n = \frac{1}{\Delta x} \minmod\left( (\matI - \matL) \bar u^n  , (\matR - \matI) \bar u^n \right). 
\end{equation}
Then the MUSCL scheme  computes the right and left states at the cell interfaces, denoted by $u^{n\pm} = (u_{3/2}^{n\pm}, \dots , u_{N_x+1/2}^{n\pm})$. They are  defined by \eqref{eq:Full_reconstruction},
\begin{align}
    u^{n\pm} = \bar u^n \pm \frac{\Delta x}{2} (u_x)^n.
\end{align}

In the case of $m$ stochastic variables, the tensor train approximation of the solution uses the representation \eqref{eq:FullyTT_def}. The previous computations are repeated on the tensor train format. Since the operation of shifting an index to the right or to the left (given by the previously defined matrices) is applied to the spatial indices, only the spatial core $G_0$ of the representation is affected. 
For example, the right-shifted core $\tilde G_0$ is given by $\tilde G_0 = \matR G_0$, and the tensor train for the shifted solution is 
\begin{equation}
    \bar{U}^n(i+1,\ibfj) 
    = \sum_{\alpha_1 \dots \alpha_m =1 }^{r_1 \dots r_m}
    \tilde G_0(i,\alpha_1) G_1(\alpha_0,j_1,\alpha_1) G_2(\alpha_1, j_2, \alpha_2) \cdots G_m(\alpha_m, j_m).
\end{equation}
The approximate derivative in the full-TT format is obtained by shifting the cores and applying the minmod limiter to the resulting tensor trains. This also requires e.g.~cross--approximation for the limiter. 
The remaining operations are linear on the tensor train and their description is omitted. 
The MUSCL scheme for the full-TT format is  summarized in Algorithm \ref{algo:muscl_fullyTT}. Finally, expectation and variance are obtained by  integration. This is achieved  by contracting the solution tensor with the rank-one tensor of the integration weights.

\section{Numerical results} \label{sec:Numerics}
In this section, we apply the method to three configurations: 
the Riemann problem for Burgers' equation, the Sod problem and the Shu-Osher problem for the Euler equations.
In the first problem, we investigate the influence of selected parameters on accuracy and compare the results with the full-TT format. We further demonstrate that the method scales effectively to a larger number of stochastic dimensions.  
The second problem illustrates the applicability of the method to a system of equations, where a parameter study and a comparison with the full-TT format are conducted.  
Finally, the third test case demonstrates that the method is capable of handling complex wave interaction problems.

The method was implemented in python with the tntorch library \cite{usvyatsov_tntorch_2022}. 

\subsection{Stochastic Burgers' equation}
We consider  Burgers' equation with uncertain initial conditions 
\begin{align} 
    \frac{\partial u}{\partial t} + \frac{1}{2}\frac{\partial u^2}{\partial x} = 0, \label{eq:Burgers}
    \quad
    u_0(x,\bfxi) = \begin{cases}
         1 + \sum_{\ell=1}^m v_{L,\ell} \, \xi_\ell, & \quad x < 0
        \\
        - 1 + \sum_{\ell=1}^m v_{R,\ell} \, \xi_\ell, & \quad x > 0
    \end{cases}  
    .
\end{align}
The stochastic variables are assumed to be independent and identically distributed: $\xi_\ell \sim \mathcal U(0,1)$, for all $\ell \in \{1,\dots,m\}$. 
Computations are performed for $\Omega_x=(-1,1)$ up to final time $T=0.35$, and we use a CFL number $\alpha=0.45$.
Let $\bfv_L=(v_{L,1},\dots,v_{L,m})^\top$ and $\bfv_R=(v_{R,1},\dots,v_{R,m})^\top$.
For now, the number $m$ of random variables is arbitrary. The value of $m$ and the coefficients $\bfv_L, \bfv_R$ will be given in each subsection. 
We focus on the case where the solution is a shock wave, hence, we assume that for all realizations of $\bfxi=(\xi_1,\dots,\xi_m)$ it holds 
\begin{equation} 
    1 + \sum_{\ell=1}^m v_{L,\ell} \, \xi_\ell > - 1 + \sum_{\ell=1}^m v_{R,\ell} \, \xi_\ell. 
\end{equation}
Under this assumption, the exact solution to \eqref{eq:Burgers} is a shock
\begin{equation}
    u(t,x,\bfxi) = \begin{cases}
        u_L(\bfxi),  & \quad x < s(\bfxi)t
        \\
        u_R(\bfxi), & \quad x > s(\bfxi)t
    \end{cases}
    ,\quad s(\bfxi) = \frac{u_L(\bfxi)+u_R(\bfxi)}{2}.
\end{equation}
In the following we focus on the case $m=3$ and study the algorithmic performance for an increasing number of (spatial and stochastic) cells. 
Then we consider the scale-up possibility by increasing $m$. 

\subsubsection{Convergence study and influence of parameters} \label{sec:Exp1_2}
We choose $m=3$ and take the following initial conditions 
\begin{equation}
    u_0(x,\bfxi) = \begin{cases}
        1 + 0.1 \xi_1 - 0.1 \xi_3,  & \quad x < 0
        \\
        - 1 + 0.1 \xi_1 - 0.1 \xi_2, & \quad x > 0
    \end{cases}
    .
\end{equation}
Next we investigate the influence of the algorithmic parameters on the performance of the algorithm. 
More precisely, the number of stochastic and spatial cells, the TT-rank and the tolerance for TT approximations are varied. We compare  the relative $L^1$-error of the expectation, 
the relative $L^1$-error of the variance,  and the number of coefficients of  the hybrid format at final time. 
The relative $L^1$-error for the expectation is defined by 
\begin{align}
    |\mathbb E|_{L^1, \text{rel}}
    = \frac{|\mathbb E_\text{TT} - \mathbb E_\text{ex}|_{L^1}}{|\mathbb E_\text{ex}|_{L^1}}.
\end{align}
The  expectation along $x$ is approximated using the midpoint rule, and the relative error is computed by
\begin{equation}
    |\mathbb E|_{L^1} 
    = \frac{1}{N_x} \sum_{i=1}^{N_x} \mathbb E[\bar u_i^n], \qquad
    |\mathbb E_{TT} - \mathbb E_\text{ex}|_{L^1}
    = \frac{1}{N_x} \sum_{i=1}^{N_x} \left( 
        \mathbb E[U_i^n] 
        - \int_{x_{i-1/2}}^{x_{i+1/2}} \mathbb E[u(x,t)]
    \right).
\end{equation}
The relative  $L^1$-error of the variance $|\text{var}|_{L^1, \text{rel}}$ is defined analogously.
The expectation for the tensor train $\mathbb E[U_i^n]$ is computed as described in Section \ref{sec:Exp_var}.  
The expectation and the variance of the exact solution are given in Appendix \ref{sec:Exact}. 
Since their expression requires nontrivial integration, the integrals are numerically approximated using the scipy.integrate package \cite{virtanen_scipy_2020}. 

The number of cells in the spatial direction is denoted by $N_x$. 
Each stochastic direction has the same discretization, and we denote by $N_\xi$ the number of cells for each stochastic direction.  
For example, in the case $N_x = N_\xi = 160$, the spatial-stochastic grid is made of $160^4$ cells. 
The TT-rank $\rmax$ corresponds to the maximal rank allowed at each rounding operation.
The tolerance $\varepsilon$ is a parameter related to the cross approximation algorithm and the rounding operations. 
It should be noted that this tolerance does not give an indication on the  global error of the tensor approximation, and that the influence of $\varepsilon$ on the total error of the algorithm is not known.
Table \ref{tab:Exp1_2} shows the result of the parameter study. The maximal rank actually observed for each simulation is denoted by $r$ in the table.

For the expectation, except for the case $\varepsilon=0.1$, all combinations of $\varepsilon$ and $\rmax$ result in  similar convergence results.
If $\varepsilon<0.1$, the order of convergence  is close to one for the expectation and close to $\frac12$  for the variance. 
Those orders are in agreement with the expected theoretical result  \cite{jin_stochastic_2017}.
For the case $\varepsilon=0.1$, the TT approximation is not accurate enough,  i.e.,the expected convergence rate is not reached. 
This behavior has also been  observed in \cite{walton_tensor-train_2025} for the full-TT discretization.
The lowest rank $\rmax=1$ seems to be as good as any larger rank when approximating the expectation. In Table \ref{tab:Exp1_2}, we see that the maximal rank $\rmax$ is reached when $N_x$ is also large and $\varepsilon$ is small. 
Finally, we note that during the simulation  we observed that for the given thresholds  the cross-interpolation did not necessarily converge. We still observe that the final errors are sufficiently small. 

\begin{figure}
\centering
\includegraphics[width=\textwidth]{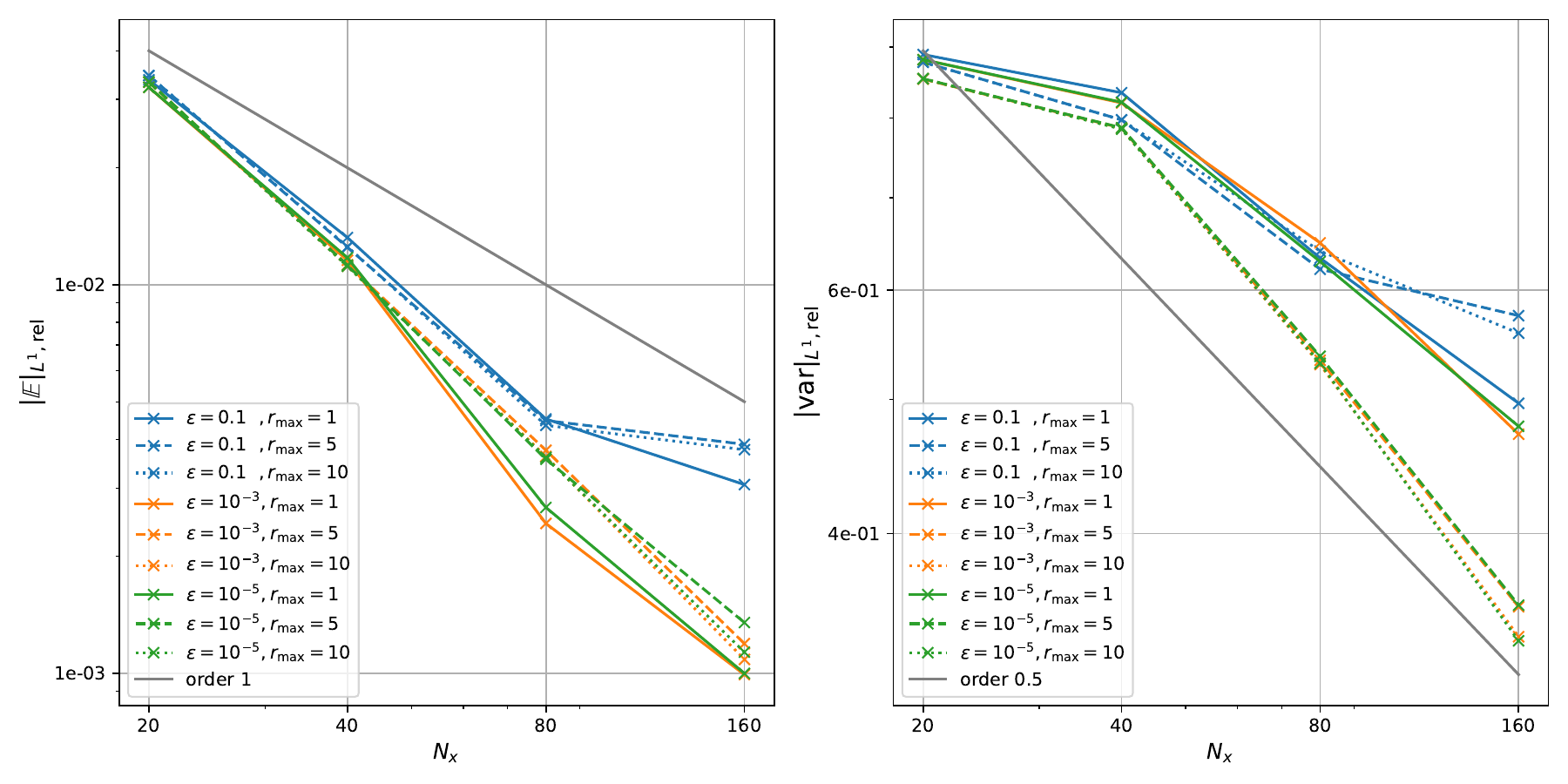}
\caption{ Relative $L^1$-errors of the expectation and variance for the shock wave computed with $N_\xi=20$ using the hybrid method.}
\label{fig:Exp2}
\end{figure}

The convergence study was then conducted for varying $N_x$ and fixed $N_\xi=20$. 
Figure \ref{fig:Exp2} shows the resulting convergence plots. 
They are very similar to the case $N_\xi=N_x$, indicating that the convergence rate is mostly driven by the choice of the accuracy in the spatial dimension $x$.

\begin{table} \label{tab:Exp1_2}
    \centering
\begin{tabular}{ | c | c | c  | c || c | c | c | c || c | c | c | c |}
    \hline   
    \multicolumn{4}{|c||}{$\varepsilon=10^{-1}$} &  \multicolumn{4}{c||}{$\varepsilon=10^{-3}$} & \multicolumn{4}{c|}{$\varepsilon=10^{-5}$} \\
    \hline
    $N_x$ & $\rmax$ & $r$ (h) & $r$ (f) & $N_x$ & $r_\text{max}$ & $r$ (h) & $r$ (f)  & $N_x$ & $r_\text{max}$ & $r$ (h) & $r$ (f) \\
\hline 
     &  1  &  1   & 1 &    &  1  &  1   & 1 &    &  1  &  1   & 1  \\ 
  20 &  5  &  2   & 1 & 20 &  5  &  2   & 3 & 20 &  5  &  5   & 5  \\ 
     &  10 &  2   & 1 &    &  10 &  2   & 3 &    &  10 &  8   & 10 \\ 
     &  30 &  2   & 1 &    &  30 &  2   & 3 &    &  30 &  12  & 13  \\ 
  \hline
     &  1  &  1   & 1 &    &  1  &  1   & 1 &    &  1  &  1   & 1  \\ 
  40 &  5  &  2   & 1 & 40 &  5  &  5   & 4 & 40 &  5  &  5   & 5  \\ 
     &  10 &  2   & 1 &    &  10 &  6   & 4 &    &  10 &  10  & 10 \\ 
     &  30 &  2   & 1 &    &  30 &  8   & 5 &    &  30 &  30  & 30  \\ 
  \hline 
     &  1  &  1   & 1 &     &  1  &  1   & 1 &    &  1  &  1   & 1  \\ 
  80 &  5  &  2   & 1 &  80 &  5  &  5   & 5 & 80 &  5  &  5   & 5  \\ 
     &  10 &  2   & 1 &     &  10 &  9   & 6 &    &  10 &  10  & 10 \\ 
     &  30 &  2   & 1 &     &  30 &  12  & 7 &    &  30 &  30  & 30  \\ 
  \hline
      &  1  &  1  & 1 &     &  1  &  1  & 1 &     &  1  &  1  & 1 \\ 
  160 &  5  &  2  & 1 & 160 &  5  &  5  & 5 & 160 &  5  &  5  & 5 \\ 
      &  10 &  2  & 1 &     &  10 &  10 & 7 &     &  10 &  10 & 10 \\ 
      &  30 &  2  & 1 &     &  30 &  15 & 7 &     &  30 &  30 & 30 \\ 
\hline 
\end{tabular}
\caption{Parameters tested for Sections \ref{sec:Exp1_2}. 
The maximal rank observed for the hybrid and the full-TT format in the case $N_\xi=N_x$ are denoted by $r$ (h) and $r$ (f), respectively.}
\end{table}

\begin{figure}
    \centering
    \includegraphics[width=\textwidth]{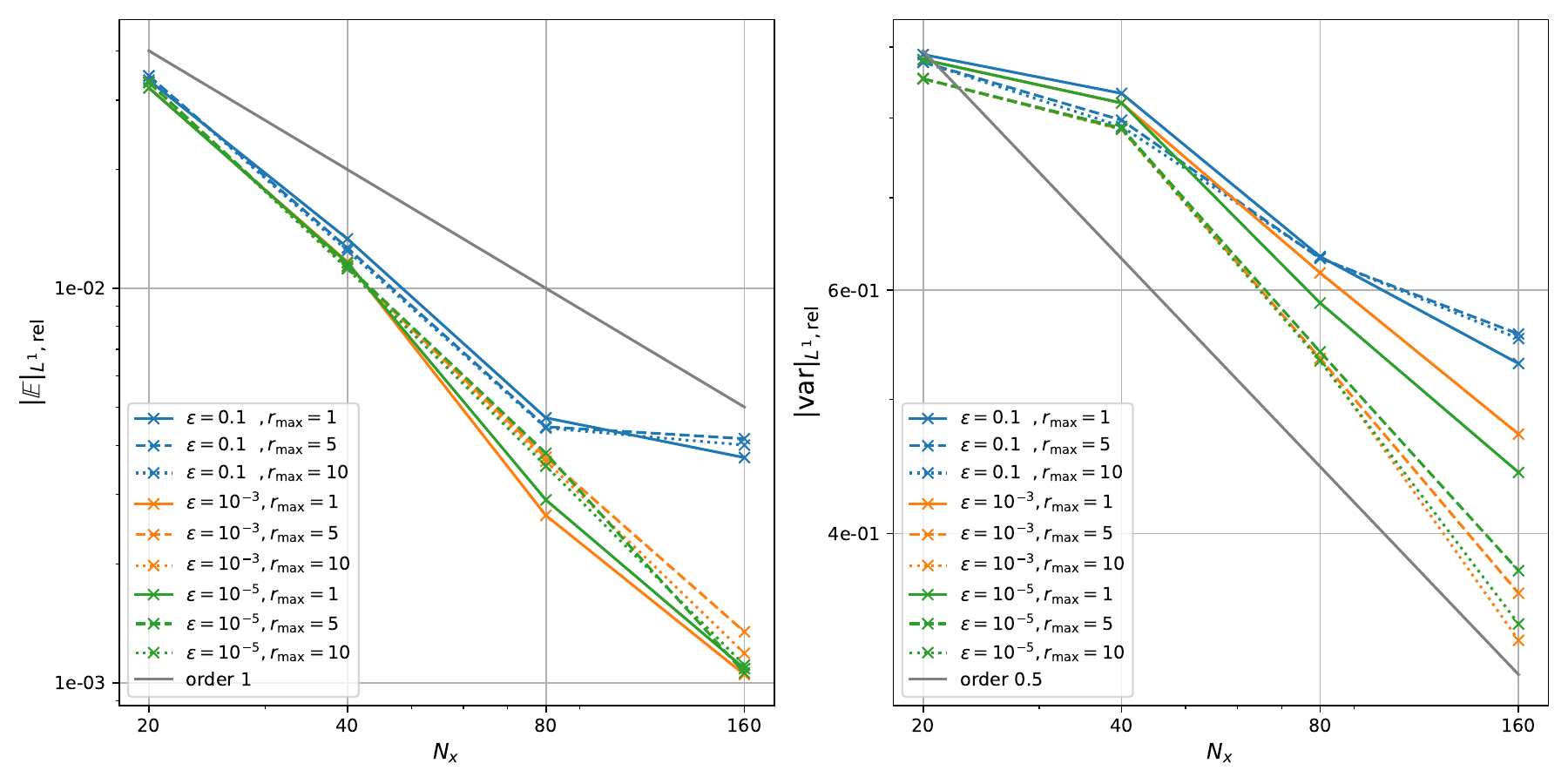}
    \caption{Relative $L^1$-errors of the expectation and variance for the shock wave computed with $N_\xi=N_x$ using the hybrid method.}
    \label{fig:Exp1}
\end{figure}

\begin{figure}
    \centering
    \includegraphics[width=\textwidth]{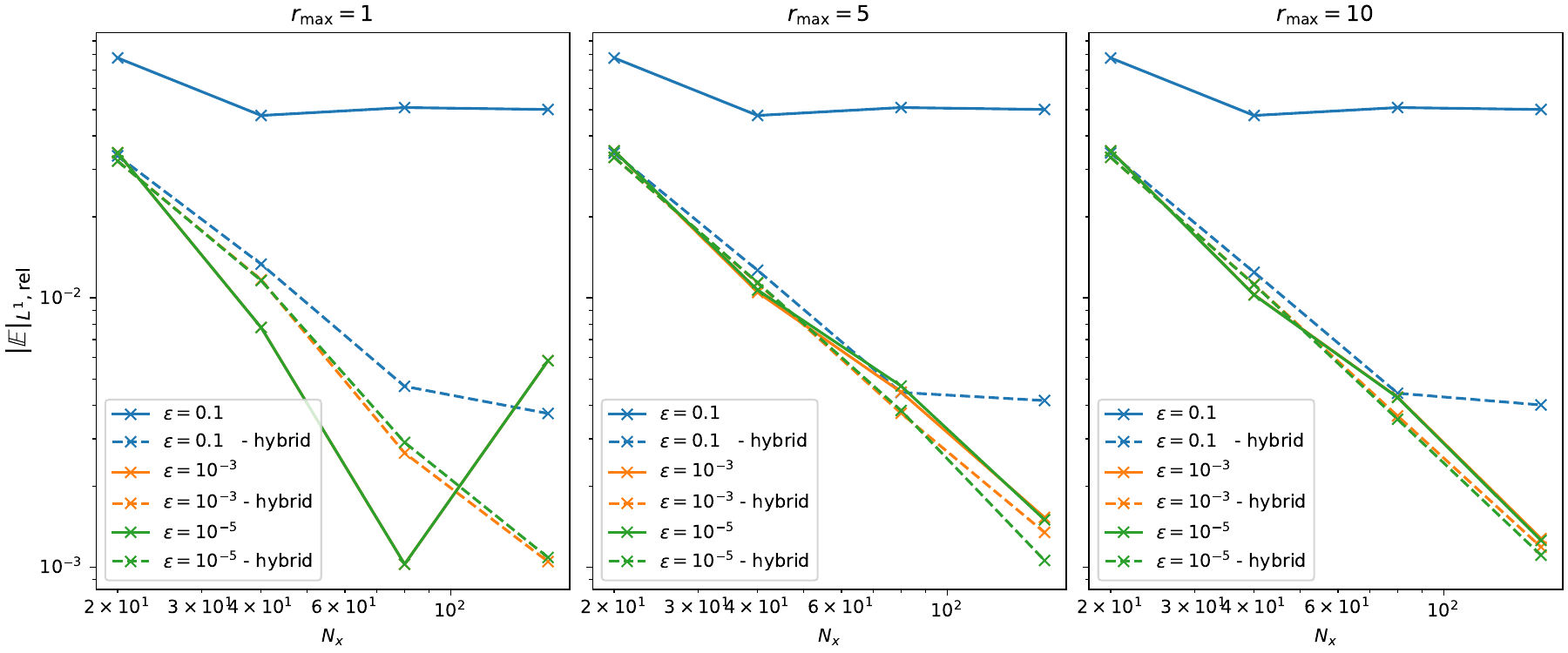}
    \caption{Relative $L^1$-error of the expectation for the shock wave using the full-TT (plain line) and hybrid (dashed line) formats.}
    \label{fig:Burgers_comparison}
\end{figure}

\begin{figure}
    \centering
    \includegraphics[width=0.5\textwidth]{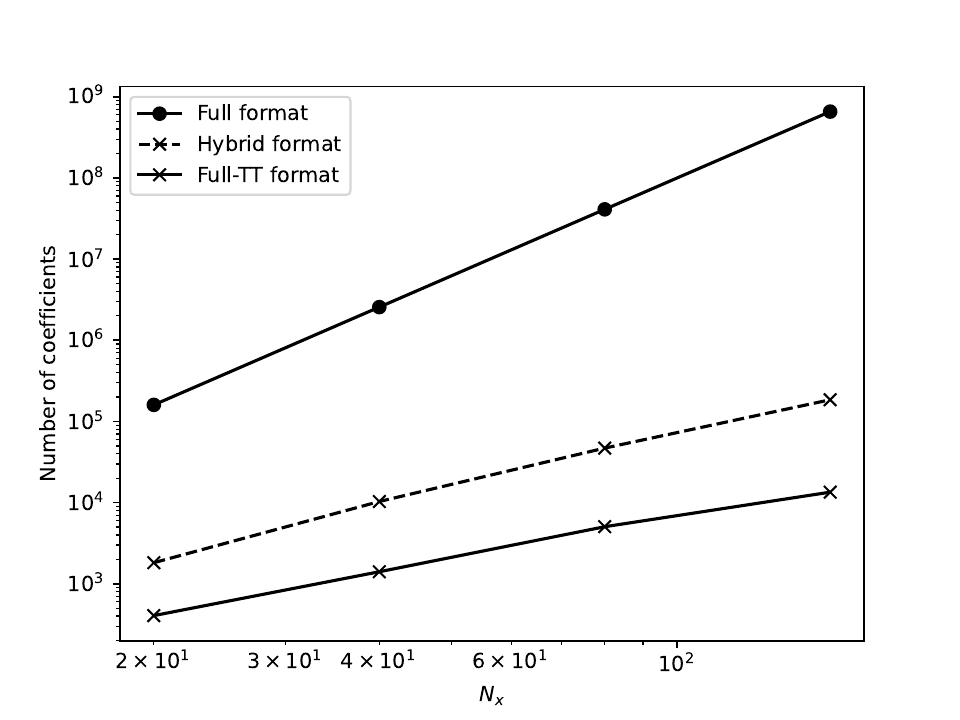}
    \caption{Number of coefficients of the solution at the final time: hybrid format, full-TT format and theoretical value for the full tensor.}
    \label{fig:Exp1_nparam}
\end{figure}

We also compare the performance of the hybrid format to the full-TT format. 
Figure \ref{fig:Burgers_comparison} shows the convergence curve of the hybrid and the full-TT format using the values of $\varepsilon$ and $\rmax$ given in Table \ref{tab:Exp1_2}.
The case $\varepsilon=0.1$ yields better results for the hybrid format than for the full-TT format.
For smaller $\varepsilon$, both methods yield comparable results. 
In this case, the full-TT method is more efficient since it requires fewer tensor coefficients: for the same maximal rank $r$, the hybrid method requires $N_x$ times more parameters. 
Figure \ref{fig:Exp1_nparam} shows how the number of coefficients at final time increases with $N_x$ for each method, compared to the theoretical number of parameters for the full discretization. 
We use $\varepsilon=10^{-3}, \rmax=5$ and $N_\xi=N_x$.
As expected, the full-TT format has a lower number of coefficients than the hybrid format, and both methods provide significant compression compared to the full tensor.

\subsubsection{Scalability results} 
\begin{figure}
    \centering
    \includegraphics[width=0.5\textwidth]{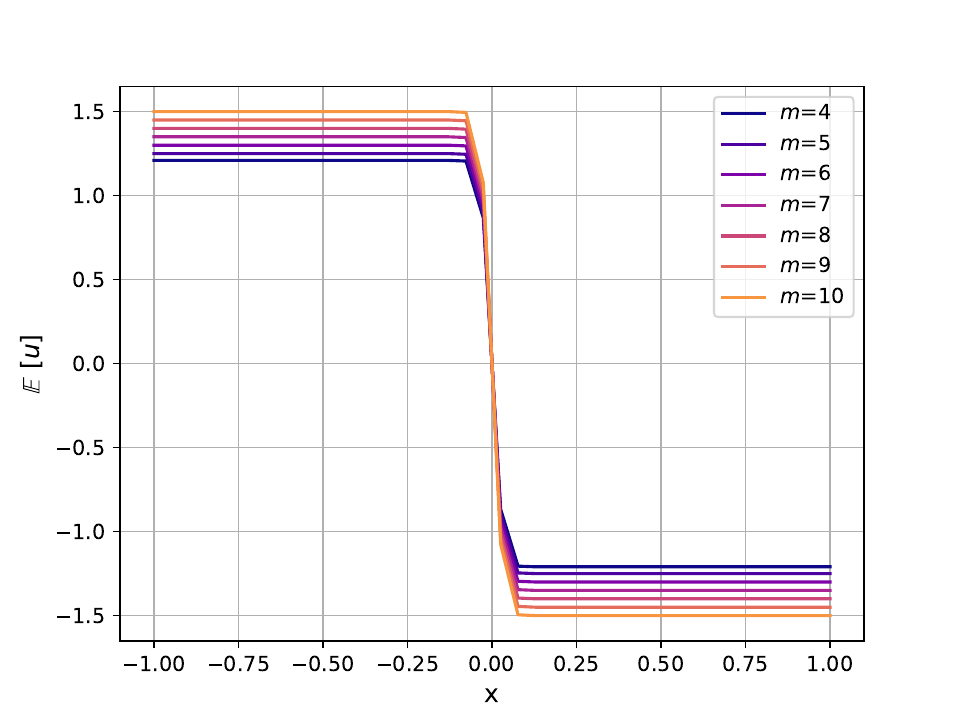}
    \caption{Expectation of the solution as a function of $x$ for various numbers of stochastic dimensions $m$.}
    \label{fig:Exp3_E}
\end{figure}

The scalability of the method is investigated for the following set of parameters: 
we fix $N_x=40,N_\xi=20, \rmax=5, \varepsilon=0.01$ and vary the number of stochastic dimensions $m$.
The coefficients $\bfv_R, \bfv_L$ are $v_{R,\ell}=0.1,~v_{L,\ell}=-0.1$ for $\ell\in\{1,\dots,m\}$. 
This ensures that the solution is a shock wave, independent of  the number of stochastic dimensions. 
Thus, the expectation will exhibit a discontinuity at $x=0$, see Figure \ref{fig:Exp3_E}.

\begin{figure}
    \centering
    \includegraphics[width=0.5\textwidth]{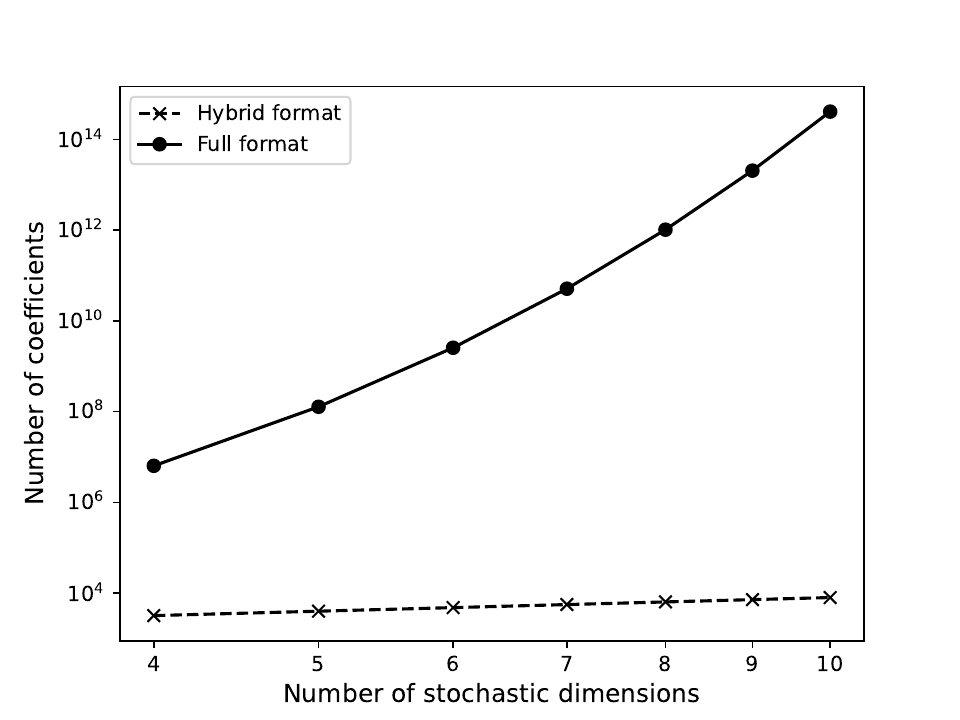}
    \caption{Number of coefficients of the solution at the final time for increasing stochastic dimensions: hybrid format and theoretical value for the full tensor}
    \label{fig:Exp3}
\end{figure}
Figure \ref{fig:Exp3} displays the number of entries in the solution at the final time, compared with the number of entries in the full tensor.  
For the hybrid format, the number of entries grows polynomially with respect to the number of stochastic dimensions, whereas for the full format it increases exponentially.
During the experiments, we observed that the cross-approximation algorithm did not always converge, with failures occurring more frequently as the number of dimensions increased.  
To assess the impact of these convergence issues on the solution, Figure \ref{fig:Exp3_E} shows the expectation of the solution. The numerical results indicate that the failed convergence does not significantly affect the solution: a shock remains clearly visible at the expected location.

\subsection{Sod problem}\label{sec:Euler}
The second experiment considers the Sod shock tube problem for the Euler equations.  
This example demonstrates that (i) the algorithm can be extended to systems of conservation laws, and (ii) the approximated solution is capable of capturing characteristic features of nonlinear hyperbolic conservation laws, such as shocks and contact waves.  
In addition, we investigate the influence of both the mesh size and the tensor rank on the accuracy of the solution.

\subsubsection{Problem description}
We consider the Euler equations in one spatial dimension,
\begin{equation} \label{eq:Euler1}
    \frac{\partial}{\partial t}
    \begin{pmatrix}
        \rho \\ \rho u \\ \rho E
    \end{pmatrix}
    +
    \frac{\partial}{\partial x}
    \begin{pmatrix}
        \rho u \\ \rho u^2 + p \\ (\rho E+p)u
    \end{pmatrix}
    = 0.
\end{equation}
Here $\rho$ is the density, $u$ the velocity, $E$ the total energy, $p$ the pressure. 
For a perfect gas with specific heat ratio $\gamma=1.4$, the pressure is given by 
\begin{equation} \label{eq:Euler2}
    p = (\gamma-1) \left( \rho E-\frac{\rho u^2}{2} \right).
\end{equation}
Note that unlike Burgers' equation, the flux of the Euler equations consists of non-polynomial functions. Hence, cross approximation is required for the  evaluation of the flux function. 
We consider the Sod shock tube problem for $x \in (0,1)$ with free boundary conditions and with $m=3$ stochastic directions, 
\begin{equation}
    \begin{pmatrix}
        \rho(0,x,\bfxi), 
        \\
        u (0,x,\bfxi), 
        \\
        p(0,x,\bfxi)
    \end{pmatrix}
    = 
    \begin{cases}
        \begin{pmatrix}
            1 + 0.1 \, \xi_1 + 0.1 \, \xi_2 + 0.05 \, \xi_3
            \\ 
            \quad -0.01 \, \xi_1 + 0.05 \, \xi_2 + 0.01 \, \xi_3
            \\ 
            1 + 0.1 \, \xi_1 - 0.01 \, \xi_2 + 0.01 \, \xi_3
        \end{pmatrix} , & 0<x<0.5
        \\[20pt]
        \begin{pmatrix}
            0.125 + 0.05 \, \xi_1 - 0.05 \, \xi_2 + 0.01 \, \xi_3
            \\ 
            0.05 \, \xi_1 - 0.01 \, \xi_2
            \\ 
            0.1 + 0.01 \, \xi_1 + 0.05 \, \xi_2 - 0.01 \, \xi_3
        \end{pmatrix} , & 0.5<x<1
    \end{cases}
    .
\end{equation}
The random variables are assumed to be independent and $\xi_\ell \sim \mathcal U(0,1)$ for all $\ell \in \{1,\dots,m\}$.  Here,  the time step is updated at each iteration such that the CFL condition holds.
For a hyperbolic system with eigenvalues $\lambda_1,\dots,\lambda_p$, the CFL condition is given by 
\begin{equation} 
    \frac{\Delta t}{\Delta x}
    \max_p |\lambda_p| < \alpha,
\end{equation}
and we use the CFL number $\alpha=0.4$. The eigenvalues of the Euler equations are 
\begin{equation}
    \lambda_1 = u-c, \ \lambda_2 = u, \ \lambda_3 = u+c,
    \quad c = \sqrt{\frac{\gamma p}{\rho}}. 
\end{equation}

\subsubsection{Parameter study}
We investigate the influence of $N_x$, $N_\xi$ and $\rmax$ on the solution.
The tolerance for TT approximations, that has been varied in the Burgers case, is now fixed at $\varepsilon=0.01$.
Since there is no analytical solution known for the stochastic Sod problem, we focus  on  qualitative results.
The expectation and variance of the primitive variables at final time $T=0.2$ are reported.
For comparison, the expectation and variance are also computed using a Monte-Carlo (MC) method and the library MultiWave 
\cite{GerhardIaconoMayMuellerSchaefer:2015}. The later implements an adaptive discontinuous Galerkin solver with $L=6$ refinement levels corresponding 
to a uniformly refined grid with $192$ cells and polynomial elements of degree 2, i.e., third-order scheme; the expectation and variance are computed using  $125,000$ MC samples.
Details can be found in \cite{herty_multiresolution_2024,pares_higher-dimensional_2024,Kolb:2024} where similar Monte-Carlo simulations have been performed.

\begin{figure}
    \centering
    \includegraphics[width=\textwidth]{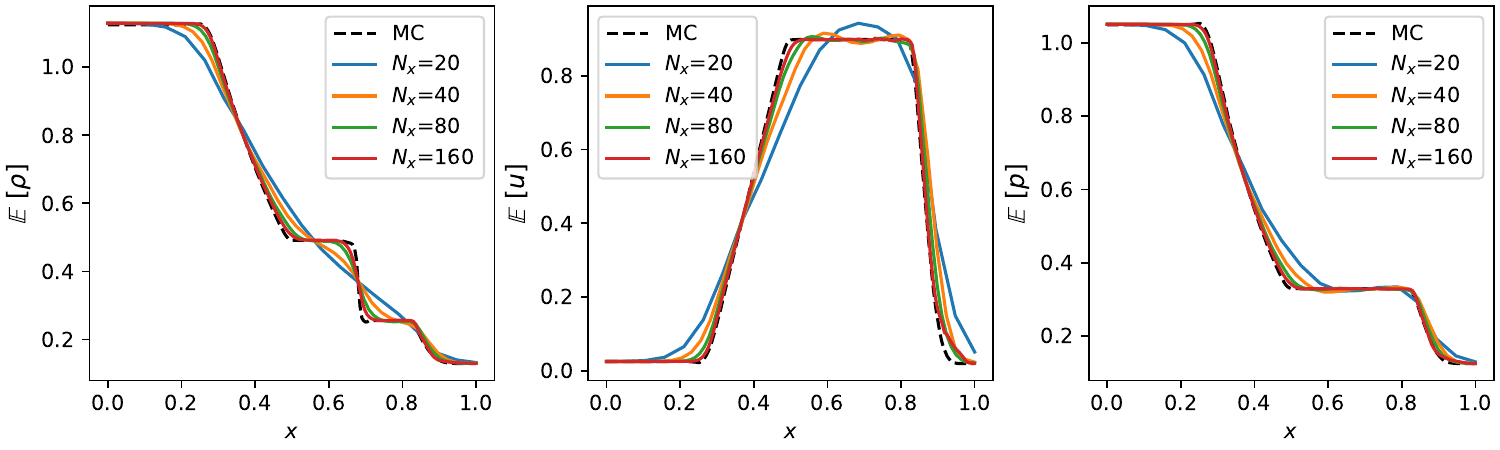}
    \caption{Expectation of $\rho, u,p$ for $\rmax=5, N_\xi=20$ and $N_x\in \{20,40,80,160\}$. The MC result (dashed line) is plotted for reference.}
    \label{fig:Euler_Nx_E}
\end{figure}

\begin{figure}
    \centering
    \includegraphics[width=\textwidth]{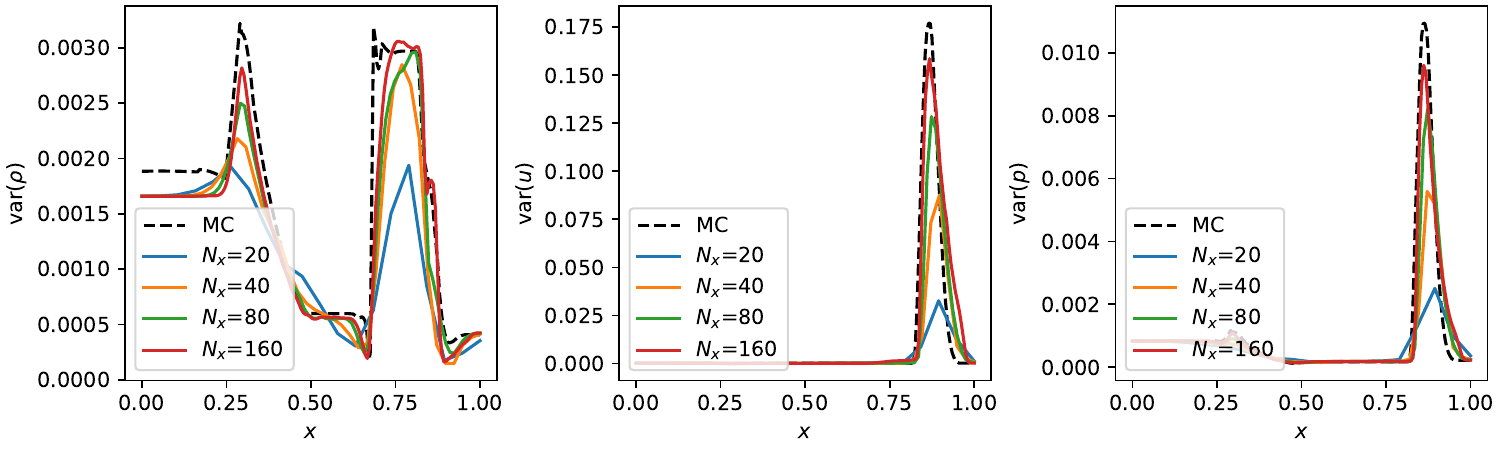}
    \caption{Variance of $\rho, u, p$ for $\rmax=5, N_\xi=20$ and $N_x\in \{20,40,80,160\}$. The MC result (dashed line) is plotted for reference}
    \label{fig:Euler_Nx_var}
\end{figure}

Figures \ref{fig:Euler_Nx_E} and \ref{fig:Euler_Nx_var}  show the expectation and variance of the primitive variable
for $r=5,N_\xi=20$ and for various mesh refinement in the physical space $N_x$, respectively. We observe convergence of the scheme when looking at the expectation 
for  $N_x=80$ and $N_x=160$, respectively. 
The qualitative behavior is also in agreement with the results of the Monte-Carlo simulations. 

In the transition between the contact wave and the shock wave at $x \in [0.65,0.85]$, there are non-physical oscillations. 
These oscillations are also present in the deterministic Euler equations  when solving by a MUSCL scheme and a leap-frog time-stepping. 
Hence, they are probably not related to the stochastic problem, and will not be further investigated in this work.
The MC results show a better resolution; this is related to the fact that they are obtained by a third-order DG solver.


\begin{figure}
    \centering
    \includegraphics[width=\textwidth]{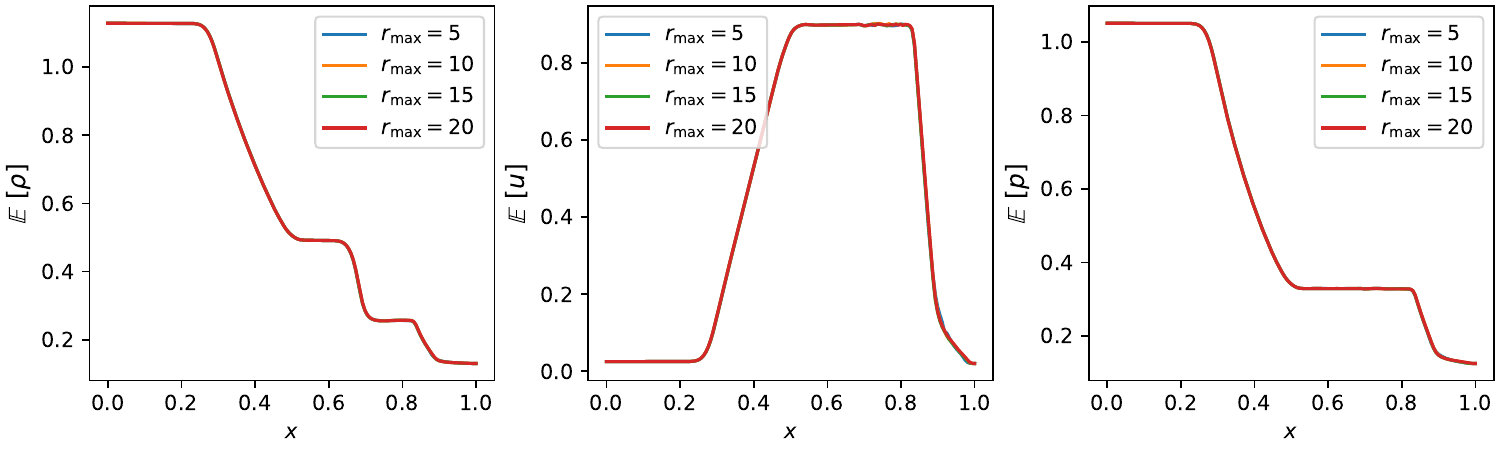}
    \caption{Expectation of $\rho, u, p$ for $N_x=160, N_\xi=20$, and $\rmax \in \{5,10,15,20\}$.}
    \label{fig:Euler_r}
\end{figure}

\begin{figure}
    \centering
    \includegraphics[width=\textwidth]{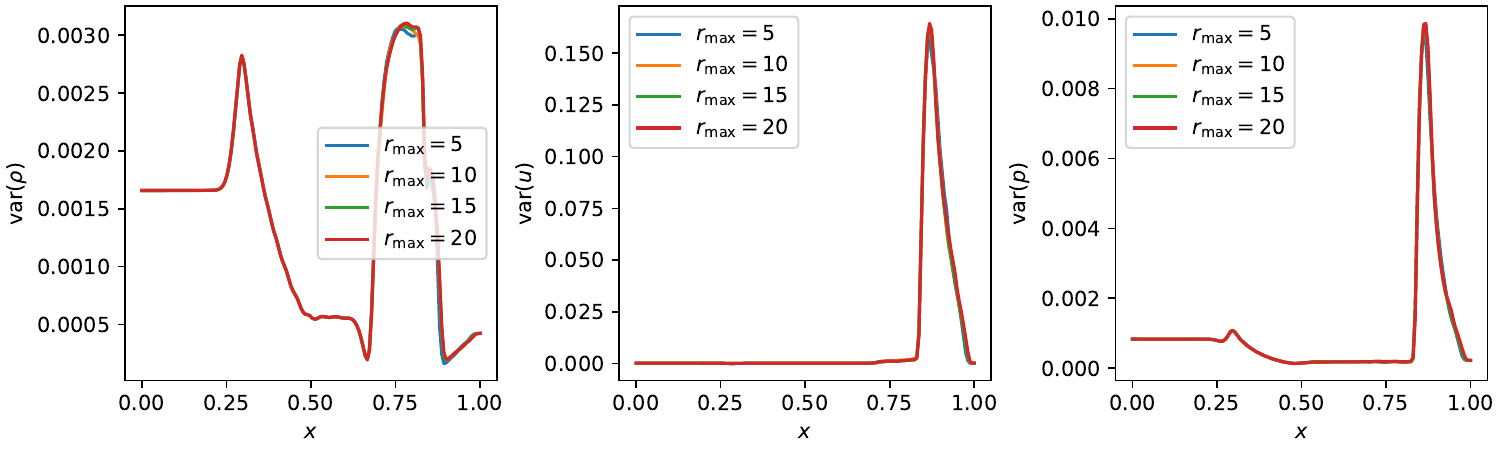}
    \caption{Variance of $\rho, u, p$ for $N_x=160, N_\xi=20$, and $\rmax \in \{5,10,15,20\}$.}
    \label{fig:Euler_r_var}
\end{figure}

\begin{figure}
    \centering
    \includegraphics[width=\textwidth]{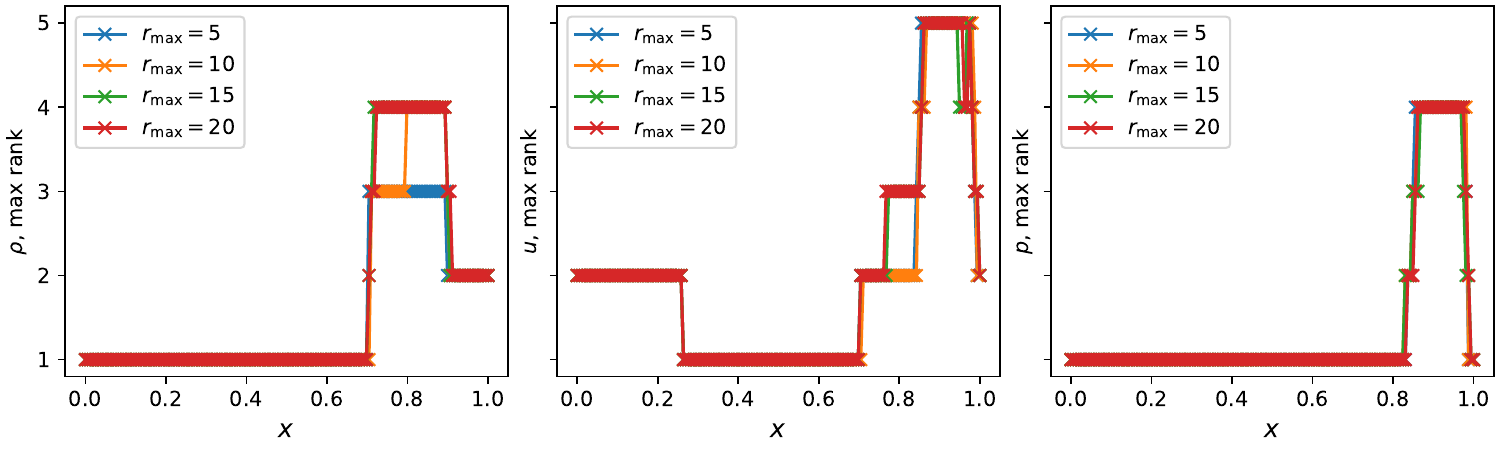}
    \caption{Maximal rank observed at final time for each spatial cell in the hybrid format of $\rho$, $u$ and $p$ for $N_x=160,N_\xi=20,\rmax\in\{5,10,15,20\}$.}
    \label{fig:Euler_rmax}
\end{figure}

To study the influence of the rank, we set $N_x=160, N_\xi=20$ and vary the maximal rank $\rmax$. The expectation and variance of $\rho$, $u$ and $p$ are shown in Figure \ref{fig:Euler_r} and in Figure \ref{fig:Euler_r_var} respectively. 
For $\rmax$ larger than 5, the results are very similar. Ranks lower than 5 lead to non-physical solutions.  
Figure \ref{fig:Euler_rmax} shows the actual maximal rank  in each spatial cell for the case $N_x=160, N_\xi=20$, $\rmax=10$ and at final time. 
In most cells, the ranks are very low -- less than three -- and  large ranks concentrate in few cells at the shock position. 
Hence, large ranks are needed only in localized regions. 


\begin{figure}
    \centering
    \includegraphics[width=\textwidth]{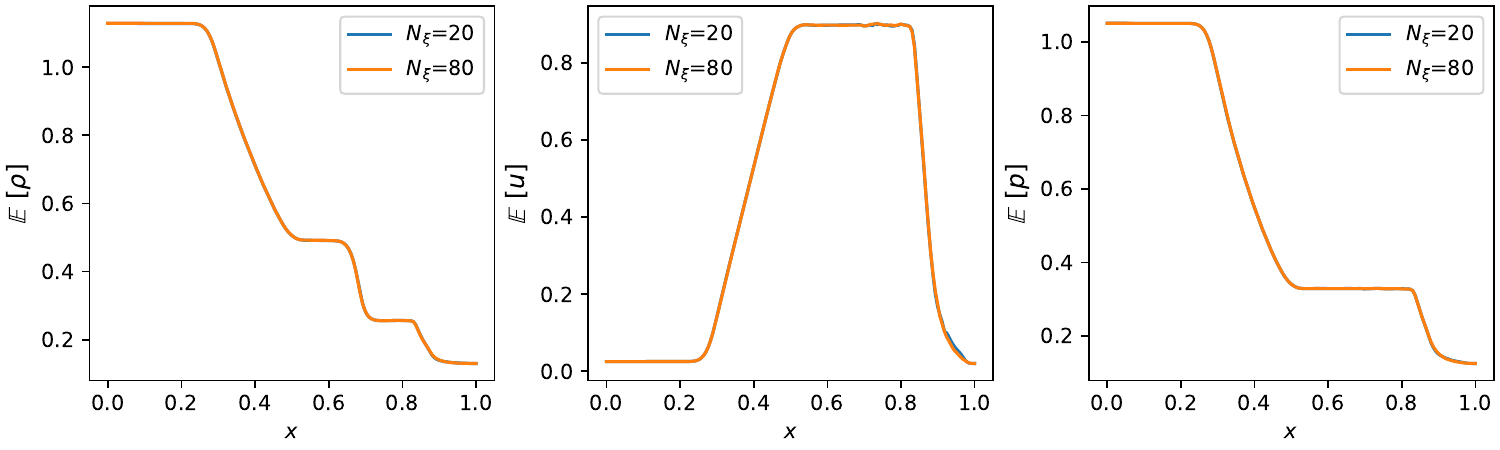}
    \caption{Expectation of $\rho, u, p$ for $N_x=160, \rmax=5$ and $N_\xi\in \{20,80\}$.}
    \label{fig:Euler_Nxi}
\end{figure}

Finally, we fix $N_x=160,r=5$ and vary $N_\xi$. Figure \ref{fig:Euler_Nxi} shows the expectation of $\rho,u,p$.
The number of stochastic cells does not seem to have a significant influence. This is probably due to the choice of the  uniform probability distribution for each random variable. 
It should be noted that the case $N_\xi=40$ yields non-physical solutions and those are not included here.

\subsubsection{Comparison with the full-TT format}
\begin{figure}
    \centering
    \includegraphics[width=\textwidth]{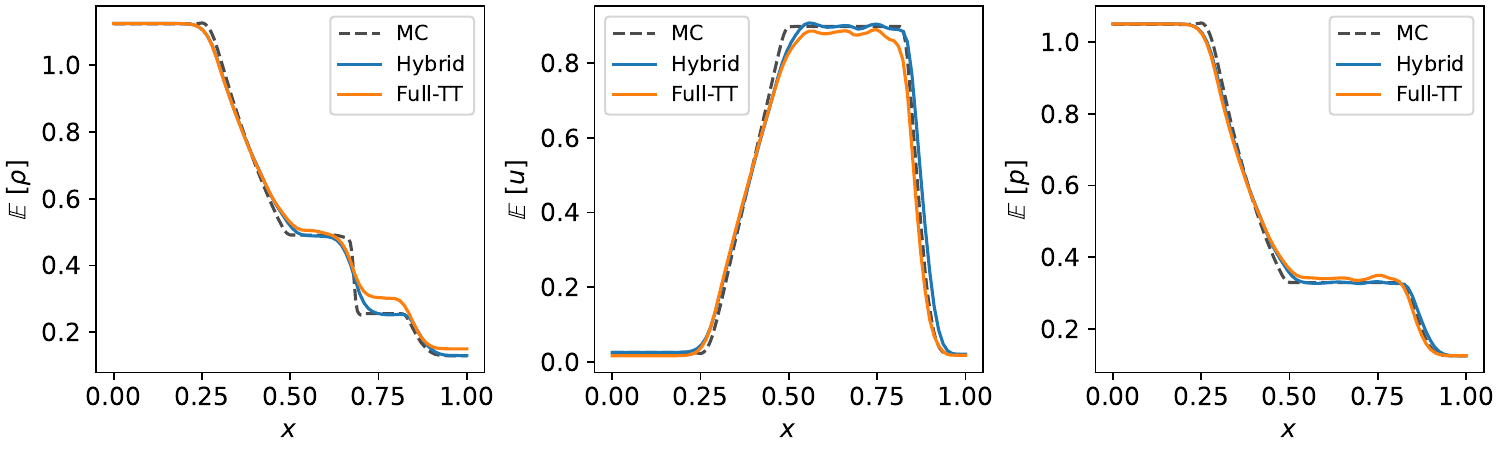}
    \caption{Expectation of $\rho$, $u$ and $p$ for the full-TT format and the hybrid format for $N_x=80, N_\xi=20$ and $\rmax=10$. The MC result (dashed line) is plotted for reference.}
    \label{fig:Euler_comparison}
\end{figure}

The results obtained with the hybrid format are compared with the full-TT formats. 
Figure \ref{fig:Euler_comparison} shows the expectation of $\rho$, $u$ and $p$ using $N_x=80, N_\xi=20,\varepsilon=10^{-3},\rmax=10$ for both formats. 
The MC result is used as a reference solution. 
We see that the hybrid format is slightly more accurate than the full-TT format, in particular for the density and around the shock.  
However, with our current implementation, the hybrid format requires more computational time than the full-TT format: it took approximately two hours to complete, against ten minutes for the full-TT format.
The number of coefficients at the final time is also much larger for the hybrid format. 
Those two aspects could be greatly improved by improving the finite volume scheme, using e.g.~a higher-order method and a non-uniform mesh in the space direction.

For the Sod problem a low rank representation is sufficient for most of the spatial domain, but higher ranks are necessary in some  parts. 
This example illustrates that the hybrid format also can be used to gain information for the  grid refinement. 

\subsection{Shu-Osher problem}
The last experiment is the Shu-Osher problem for the Euler system \eqref{eq:Euler1}-\eqref{eq:Euler2}. 
In this test case, a shock wave interacts with a sine wave in the density. The solution includes a shock wave and small-scale oscillations that need to be resolved by the scheme. The classical problem is extended by two random variables such that for a realization  $\xi_1=\xi_2=0$ we recover the original problem. The initial condition for the random problem in  primitive variables are given by
\begin{equation}
    \begin{pmatrix}
        \rho(0,x,\bfxi), 
        \\
        u (0,x,\bfxi), 
        \\
        p(0,x,\bfxi)
    \end{pmatrix}
    = 
    \begin{cases}
        \begin{pmatrix}
            3.857143 + 0.1 \xi_1
            \\ 
            2.629369 + 0.1 \xi_1
            \\ 
            10.33333 + \xi_1
        \end{pmatrix} , & 0<x<\frac{1}{8}
        \\[20pt]
        \begin{pmatrix}
            1 + 0.2(1+\xi_2) \sin(16\pi x)
            \\ 
            0.0
            \\ 
            1. + 0.1 \xi_2
        \end{pmatrix} , & \frac{1}{8}<x<1
    \end{cases}
    .
\end{equation}

The simulation is run until time $t=0.13$. The following parameters are used for the discretization: $N_x=1600, N_\xi=20,\rmax=5, \varepsilon=0.01$. 
We also run tests  with a smaller maximal rank  $\rmax$ and a larger $\varepsilon$, but the solution yields (non-physical) negative densities. Therefore, those results are omitted. The relatively high rank $\rmax$ indicates that the complex wave pattern does not allow for a strong data compression of the solution. 
\par
\begin{figure}
    \centering
    \includegraphics[width=\textwidth]{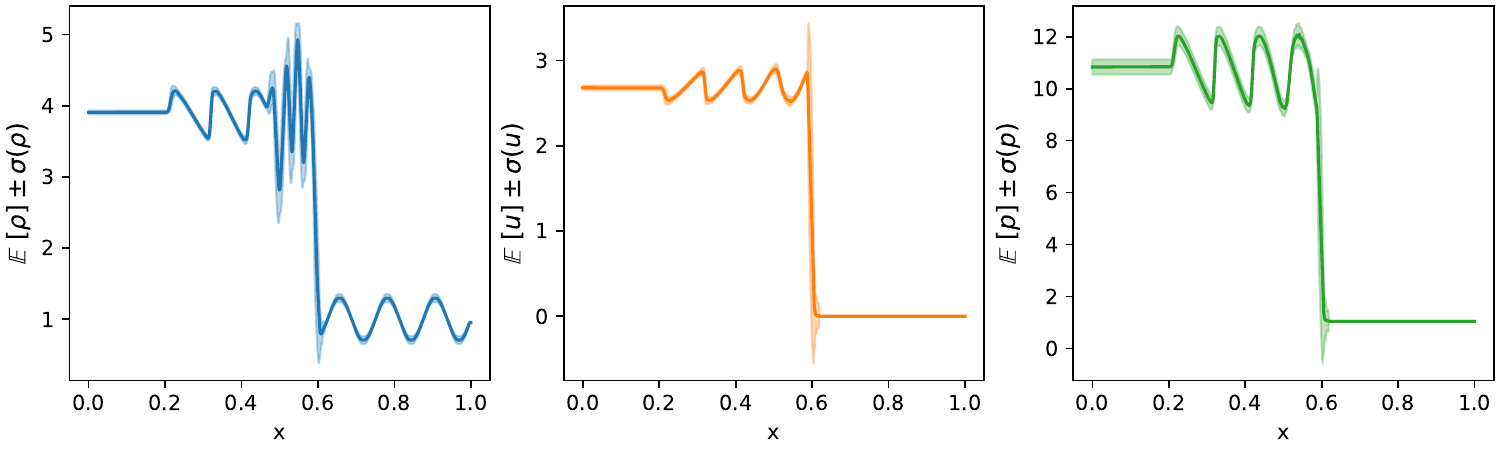}
    \caption{Expectation and standard deviation for $\rho$, $u$ and $p$ at final time.
    The results were obtained with $N_x=1600, N_\xi=20,\rmax=5, \varepsilon=0.01$.}
    \label{fig:Euler_shu}
\end{figure}

\begin{figure}
    \centering
    \includegraphics[width=\textwidth]{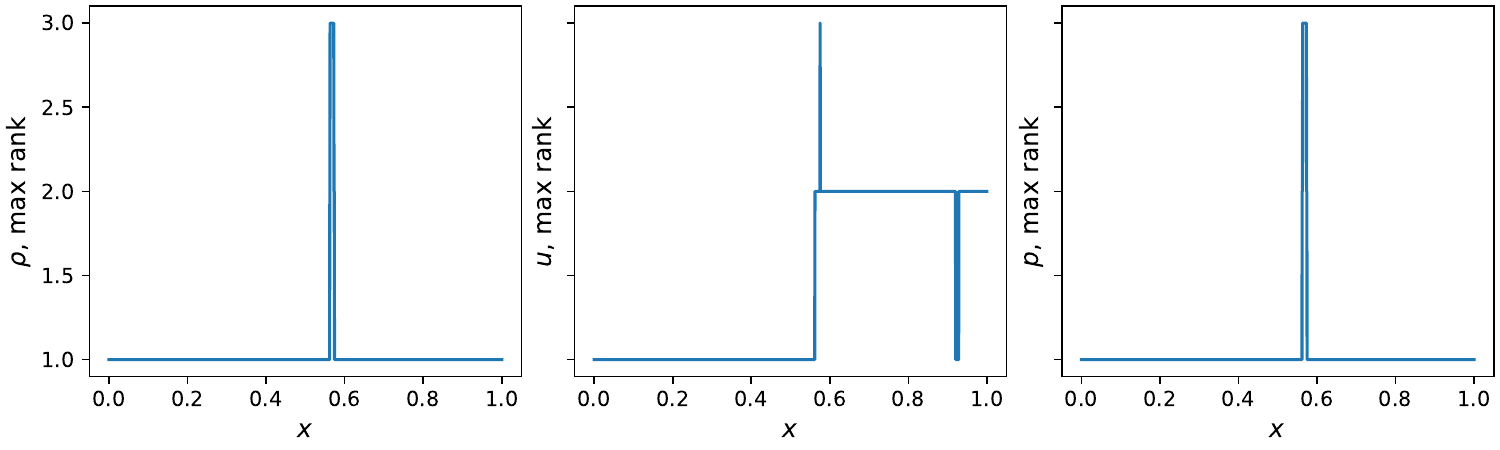}
    \caption{Maximal rank in each spatial cell for $\rho,u,p$ at final time. 
    The results were obtained with $N_x=1600, N_\xi=20,\rmax=5, \varepsilon=0.01$.}
    \label{fig:Euler_shu_rmax}
\end{figure}

 
Figure \ref{fig:Euler_shu} shows the expectation and the standard deviation $\sigma$ for $\rho$, $u$ and $p$.  We observe that the standard deviation is large near the shock wave. 
Figure \ref{fig:Euler_shu_rmax} shows the maximal ranks in each spatial cell at the final time for $\rho, u,$ and $p$, respectively.  
Interestingly, most cells have a very low rank ($r\leq2$). The largest ranks are, as in the Sod case, concentrated near the shock. This is  expected due to the low regularity in those cells. The oscillating wave patterns can also be represented by low rank structures.   The maximal rank at intermediate times are plotted in Appendix \ref{sec:A_Shu}. 
For the  velocity, the region to the right of the shock has a large rank compared with the region  to the left of the shock. As before, we observe that ranks are  low for most spatial cells and that the rank increases only near the shock wave.  However, the maximal rank is still small compared to the dimension of the problem. This  shows the usability of the hybrid format even for more complex test cases. 

\section{Conclusion and Outlook}
In this paper, we have introduced a new stochastic finite volume method capable of handling a large number of stochastic parameters.  
The numerical experiments confirm the feasibility of the approach. A comparison with the full-TT format shows that both methods produce qualitatively similar results. In the current implementation, however, the full-TT format is more efficient, both in terms of runtime and the maximal number of coefficients.
Our present implementation, based on regular grids and without parallelization, could be  enhanced. A distinctive feature of the proposed method is that the rank may vary across spatial cells, enabling adaptation to the local structure of the solution. Moreover, the approach generalizes naturally to arbitrary finite-volume schemes.

So far, numerical simulations have been carried out only in one spatial dimension.  
It is important to emphasize that, in the field of conservation laws, the complexity of the problem does not necessarily increase with the number of spatial dimensions. The principal challenges arise from the loss of regularity and the interaction of wave patterns, which are the dominant factors contributing to the complexity of such problems.

This work opens new avenues for combining the tensor train format with accurate numerical methods for hyperbolic conservation laws.  
Several aspects of the proposed method merit further investigation:
\begin{itemize}
	\item The choice of appropriate algorithm parameters -- such as the TT-rank or the tolerance for TT approximations -- is very likely to be problem-dependent. 
    Hence, the algorithm should be tested on other representative problems and for various probability densities. 
	
	\item Convergence estimates combining the SFV and the TT approximation are yet to be developed. 
    Our results indicate that the tolerance for the TT approximation has some impact on the convergence rate. Their theoretical estimates is also an ongoing work.
	
    \item In the Sod problem  and the Shu-Osher problem, it appears that the ranks are generally kept very low and that larger ranks are required near a discontinuity.
    The correlation between spatial and stochastic features is an interesting topic for future research.
\end{itemize}

The proposed algorithm is presented as a proof of concept.  
It can be readily combined with various high-performance computing techniques. In particular, the algorithm is straightforward to parallelize, since the numerical flux is computed independently within each spatial cell.
Furthermore, the observations from the Sod and Shu–Osher problems suggest a promising new direction for mesh adaptation: regions in the spatial domain where the tensor rank increases are precisely those where mesh refinement would be most beneficial. A large rank may therefore serve as an indicator for refinement.
In particular, the mesh adaptation strategies introduced in \cite{herty_multiresolution_2024, HertyKolbMueller:2024, Kolb:2024} for both physical and stochastic spaces could be adapted to the present setting.

\appendix

\section{Algorithms} \label{sec:Algo}
We describe the method in pseudo-code. 
The general algorithm for the hybrid format is presented in Algorithm \ref{algo:hybrid} and the reconstruction is presented in Algorithm \ref{algo:reconstruction}.
The general algorithm and reconstruction for the full-TT format are presented in Algorithm \ref{algo:muscl_fullyTT} and Algorithm \ref{algo:reconstruction_fullyTT}, respectively.
The CFL number is denoted by $\alpha$.

In the full-TT format, the left and right numerical fluxes are computed as follows: 
the numerical flux at each interface is first computed using the global operations, 
then the left-shifted flux is computed by applying the matrix $\matL$ to the spatial core of the flux. Those operations are described in Algorithm \ref{algo:flux_fullyTT}.

Every time cross interpolation is used, it is made visible in the pseudo-code with the syntax '$\Call{cross}{f(T1, T2, \dots)}$'. 
Here, 'f' is the non-polynomial function and 'T1, T2, ...' are its arguments in the tensor train format.

\begin{remark}
    In the Algorithms \ref{algo:hybrid} and \ref{algo:muscl_fullyTT}, the maximum of $|f'(\bar U^0)|$ can be computed directly in the tensor train format. 
    In tntorch -- the tensor train library used for this work -- the minimum of a tensor train is computed by minimizing the objective function $x \mapsto \tan(\pi/2 - x)$, where the objective function is evaluated with cross approximation. 
\end{remark}

\begin{algorithm}
\caption{MUSCL scheme and forward Euler for the hybrid format} \label{algo:hybrid}
\begin{algorithmic}
    \FORALL{$i \in \{1, \dots, N_x N_q\}$}
        \STATE $U^0_{i} \gets \Call{cross}{u_0(x_{i},\xi_1,\dots,\xi_m)}$
        \COMMENT{Evaluate $u_0$ for the midpoint rule}
    \ENDFOR
    \STATE $\bar U^0 \gets \Call{Average}{U^0}$
    \STATE $\Delta t \gets \alpha \dfrac{\Delta x}{\max |f'(\bar U^0)|}$. 
    \FORALL{$k \in 0, \dots, N_t$} 
        \FORALL{$i \in \{1, \dots, N_x\}$}
            \STATE $U^{n  \pm}_{i \pm 1/2} \gets \Call{Reconstruction}{\bar U_i^n } $ 
            \COMMENT{Reconstruct values at interfaces, see Algorithm \ref{algo:reconstruction} }
            \STATE $a \gets \Call{cross}{\max(|f'(U_{i+1/2}^+)|, |f'(U_{i+1/2}^-)|)}$
            \COMMENT{Compute the numerical velocity}
            \STATE
            \STATE $\bar H_{i\pm1/2} \gets \dfrac{f(U_{i\pm1/2}^+) + f(U_{i\pm1/2}^-)}{2 \Delta x}
                - \dfrac{a (U_{i\pm1/2}^+ - U_{i\pm1/2}^-)}{2 \Delta x}$ 
            \COMMENT{Compute left and right fluxes}
            \STATE $\bar U_i^{n +1} \gets \bar U_i^n  - \frac{\Delta t }{\Delta x}(\bar H_{i+1/2}-\bar H_{i-1/2})$
            \COMMENT{Update the solution}
        \ENDFOR
    \ENDFOR
\end{algorithmic}
\end{algorithm}

\begin{algorithm}
    \caption{Reconstruction} \label{algo:reconstruction} 
\begin{algorithmic}
    \FORALL {$i \in \{1, \dots, N_x\}$}
	\STATE $(U_x)_i^n  \gets \Call{cross}{\minmod (U_i - U_{i-1} , U_{i+1} - U_i)}$
    \STATE $U^{n +}_{i+1/2} \gets \bar U_{i+1}^n - \frac{1}{2} (U_x)_{i+1}^n$  
    \STATE $U^{n -}_{i+1/2} \gets \bar U_i^n + \frac{1}{2} (U_x)_i^n$  
	\ENDFOR
\RETURN $U^{n +}_{i+1/2}, U^{n -}_{i+1/2}$
\end{algorithmic}
\end{algorithm}

\begin{algorithm}
    \caption{MUSCL scheme in the full-TT format} \label{algo:muscl_fullyTT}
\begin{algorithmic}
    \STATE $U^0 \gets \Call{cross}{u_0(x_{i},\xi_1,\dots,\xi_m)}$
        \COMMENT{Evaluate $u_0$ for the midpoint rule}
    \STATE $\bar U^0 \gets \Call{Average}{U^0}$
    \STATE  $\Delta t \gets \alpha  \dfrac{\Delta x}{\max |f'(\bar U^0)|}$. 
    \FORALL{$k \in 0, \dots, N_t$} 
        \STATE $U^{n  \pm} \gets \Call{ReconstructionFTT}{\bar U^n }$ 
        \COMMENT{Reconstruct values at interfaces, see Algorithm \ref{algo:reconstruction_fullyTT} }
        \STATE $\bar H_L, \bar H_R \gets \Call{NumFluxFTT}{U^\pm}$ 
        \COMMENT{Compute the left and right fluxes, see Algorithm \ref{algo:flux_fullyTT}}
        \STATE $\bar U^{n +1} \gets \bar U_i^n  - \frac{\Delta t }{\Delta x}(\bar H_R-\bar H_L)$
        \COMMENT{Update the solution}
    \ENDFOR
\end{algorithmic}
\end{algorithm}

\begin{algorithm}
    \caption{Reconstruction in full-TT format} \label{algo:reconstruction_fullyTT} 
\begin{algorithmic}
    \STATE $G_0, G_1, \dots, G_m \gets \Call{cores}{\bar U^n}$ 
        \COMMENT{Collect the TT-cores of the solution}
    \STATE $\tilde G_0^L, \ \tilde G_0^R \gets (\matI - \matL) G_0, \ (\matR - \matI) G_0$
    \STATE $U^n_L, \ U^n_R \gets 
    \Call{TensorTrain}{[\tilde G_0^L, G_1, \dots, G_m]}, \ \Call{TensorTrain}{[\tilde G_0^R, G_1, \dots, G_m]}$
        \COMMENT{Create left- and right-shifted solution with the shifted cores}
    \STATE $(U_x)^n  \gets \Call{cross}{\minmod (U_L^n , U_R^n)}$
    \STATE $U^{n \pm} \gets \bar U^n \pm \frac{1}{2} (U_x)^n$   
\RETURN $U^{n +}, U^{n -}$
\end{algorithmic}
\end{algorithm}

\begin{algorithm}
    \caption{Flux computation in full-TT format} \label{algo:flux_fullyTT} 
\begin{algorithmic}
    \STATE $a \gets \Call{cross}{\max(|U^{n+}|, |U^{n-}|)}$ 
        \COMMENT{Compute the numerical velocity}
    \STATE $\bar H \gets \frac{f(U^{n+})+U^{n-}}{2} - a \frac{U^{n+}-U^{n-}}{2}$ 
    \STATE $G_0, G_1, \dots, G_m \gets \Call{cores}{\bar H}$   \COMMENT{Collect the TT-cores of the flux}
    \STATE $\tilde G_0^L \gets \matL G_0$ 
        \COMMENT{Shift the spatial core to one index to the left}
    \STATE $\bar H_L \gets \Call{TensorTrain}{[\tilde G_0^L, G_1, \dots, G_m]}$  
        \COMMENT{Create left-shifted flux from the modified cores}
\RETURN $\bar H_L, \bar H$
\end{algorithmic}
\end{algorithm}

\section{Expectation and variance for the exact solution} \label{sec:Exact}
We give here the expression for the expectation and variance of the exact solution to the Burgers' equation, for the shock case and for $m \geq 2$.  

We assume that $v_m=(v_{L,m}+v_{R,m})/2\neq0$ and define 
\begin{align}
    \hat\bfxi=(\xi_1,\dots,\xi_{m-1}), 
    \quad
    \hat{\bf v}_L=(v_{L,1},\dots,v_{L,m-1}),
    \\
    \hat{\bf v}_R=(v_{R,1},\dots,v_{R,m-1}),
    \quad
    \hat{\bf v}=(\hat{\bf v}_L+\hat{\bf v}_R)/2.
\end{align}
Let $c(t,x,\hat\bfxi)=(x/t-\hat{\bf{v}}\cdot\hat\bfxi)/v_m$ and
\begin{equation}
    a(t,x,\hat \bfxi) = \begin{cases}
        1,  & c(\hat \bfxi,x,t) > 1,
        \\
        c(t,x,\hat \bfxi),  & 0 < c(t,x,\hat \bfxi) < 1,
        \\
        0,  &  c(t,x,\hat \bfxi) < 0
    \end{cases}
    = \min\left(1, \max\left(0, c(t,x,\hat \bfxi)\right)\right).
\end{equation}

The expectation is 
\begin{multline}
    \E[u] = \int_{[0,1]^{m-1}} \int_{\xi_m = 0}^{a(\hat \bfxi)} 
          u_R(\bfxi) p_1(\xi_1) \dots p_m(\xi_m)  
    \diff \xi_m \diff \hat \bfxi 
    \\
    +  \int_{[0,1]^{m-1}} \int_{\xi_m = a(\hat \bfxi)}^{1} 
        u_L(\bfxi) p_1(\xi_1) \dots p_m(\xi_m)       
    \diff \xi_m \diff \hat \bfxi .
\end{multline}
For a uniform probability density over $[0,1]$, i.e., $p_1 = \dots = p_m = 1$ , and using the expression of $u_L,u_R$ and ${\bf v}_L, {\bf v}_R$, we have
\begin{equation}
    \E[u]
    = \int_{[0,1]^{m-1}}
          (-1+ \hat {\bf v}_R \cdot \hat \bfxi) a + v_{R,m} \frac{a^2}{2}
    \notag
        + (1 + \hat {\bf v}_L \cdot \hat \bfxi) + \frac{v_{L,m}}{2} 
        - (1 + \hat {\bf v}_L \cdot \hat \bfxi) a - v_{L,m} \frac{a^2}{2} 
    \diff \hat \bfxi .
\end{equation}
For the variance we start by computing
\begin{multline}
    \E[u^2] 
    = \int_{[0,1]^{m-1}}
    \left( 
      (-1 + \hat {\bf v}_R \cdot \hat \bfxi)^2 a
    + (-1 + \hat {\bf v}_R \cdot \hat \bfxi) v_{R,m} a^2
    + \frac{v_{R,m}^2}{3} a^3
    \right)
    \\ 
    -
    \left( 
      (1 + \hat {\bf v}_L \cdot \hat \bfxi)^2 a
    + (1 + \hat {\bf v}_L \cdot \hat \bfxi) v_{L,m} a^2
    + \frac{v_{L,m}^2}{3} a^3
    \right) 
    \\ 
    +
    \left( 
      (1 + \hat {\bf v}_L \cdot \hat \bfxi)^2 
    + (1 + \hat {\bf v}_L \cdot \hat \bfxi) v_{L,m} 
    + \frac{v_{L,m}^2}{3} 
    \right)   
    \diff \hat \bfxi.
\end{multline}
The variance is then obtained with $\var [u(t,x)] = \E[u(t,x)^2]-\E[u(t,x)]^2$.

\section{Shu-Osher problem: maximal ranks at various timesteps} \label{sec:A_Shu}
\begin{figure}
    \centering
    \includegraphics[width=\textwidth]{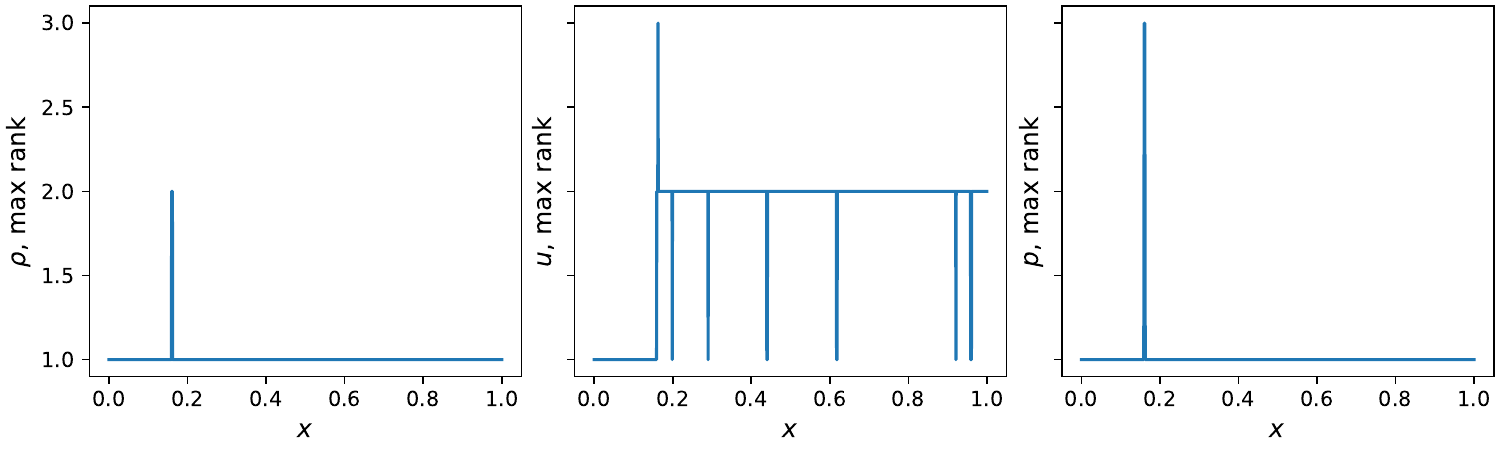}
    \caption{Maximal rank in each spatial cell for $\rho$, $u$ and $p$ at initial time. }
    \label{fig:Euler_shu_rmax0}
\end{figure}

\begin{figure}
    \centering
    \includegraphics[width=\textwidth]{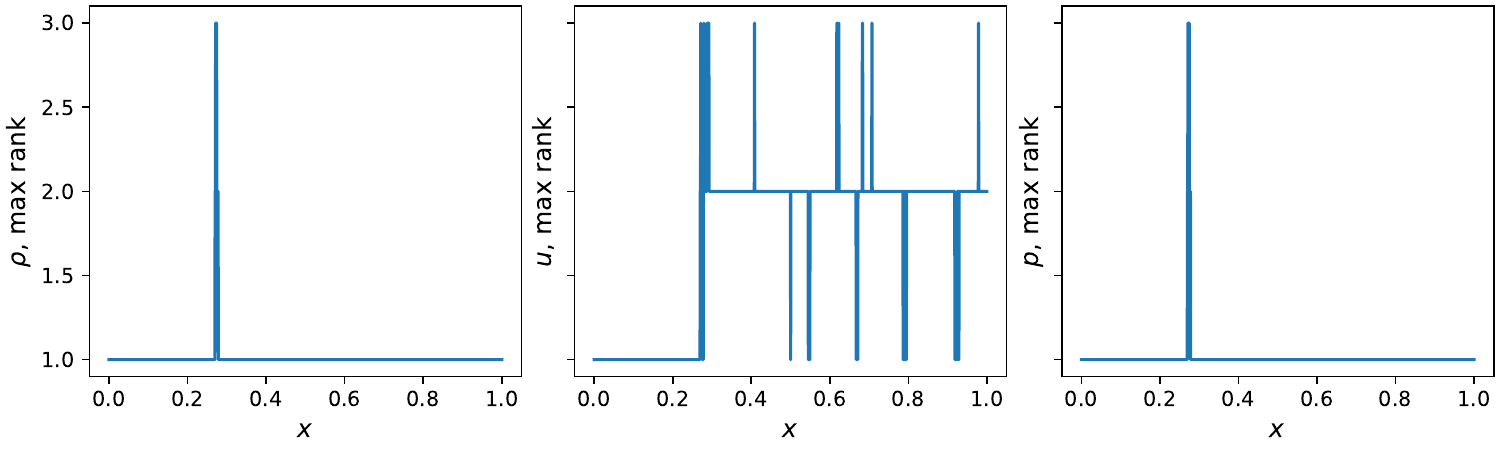}
    \caption{Maximal rank in each spatial cell for $\rho$, $u$ and $p$ at time $t=0.04$.}
    \label{fig:Euler_shu_rmax3}
\end{figure}

\begin{figure}
    \centering
    \includegraphics[width=\textwidth]{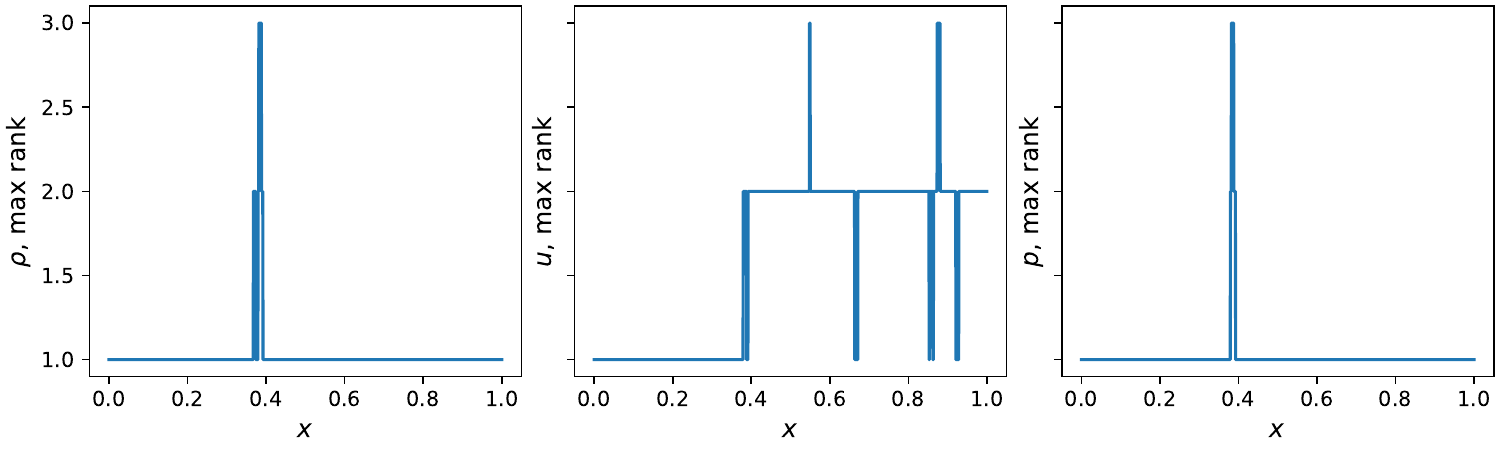}
    \caption{Maximal rank in each spatial cell for $\rho$, $u$ and $p$ at time $t=0.07$.}
    \label{fig:Euler_shu_rmax6}
\end{figure}

\begin{figure}
    \centering
    \includegraphics[width=\textwidth]{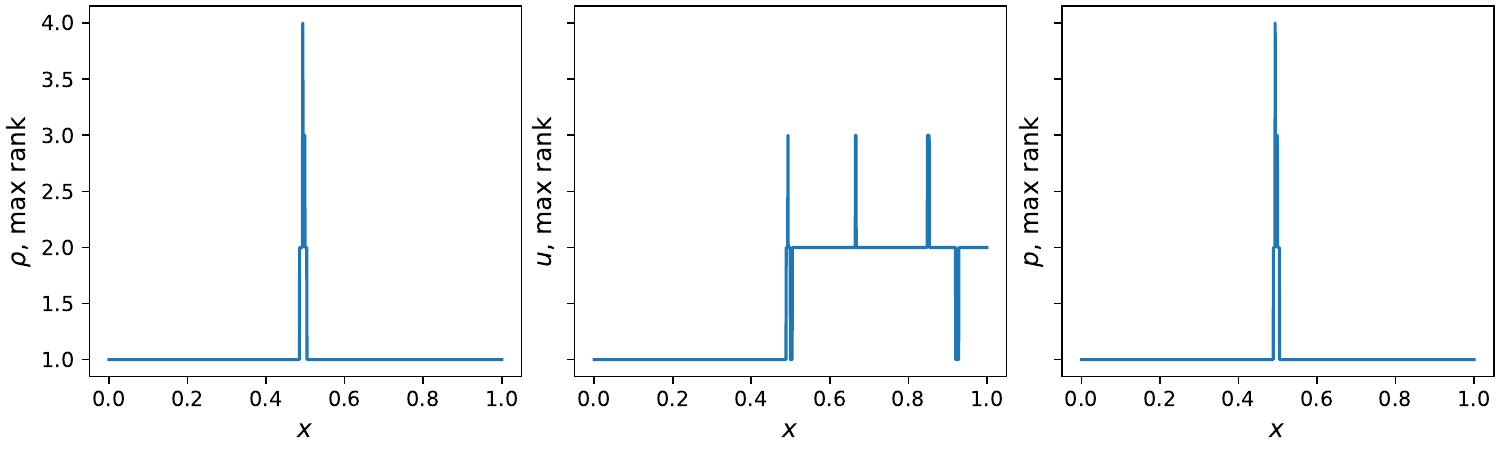}
    \caption{Maximal rank in each spatial cell for $\rho$, $u$ and $p$ at time $t=0.1$.}
    \label{fig:Euler_shu_rmax9}
\end{figure}

Figures \ref{fig:Euler_shu_rmax0}–\ref{fig:Euler_shu_rmax9} show the maximal rank in each spatial cell for the solution $\rho$, $u$ and $p$ of the Shu-Osher problem at 
times $t \in \{0.0, 0.04, 0.07, 0.1\}$. The results were obtained with $N_x=1600, N_\xi=20,\rmax=5$ and $\varepsilon=0.01$.

\section*{Acknowledgments}
We would like to thank Adrian Kolb for providing the Monte-Carlo results, and the anonymous reviewer for his useful comments.

\printbibliography 

\end{document}